\documentclass[review]{elsarticle}
  \usepackage{amsfonts, amssymb, graphicx, a4, epic, eepic, epsfig,
  setspace, lscape,color}
  \usepackage{amsmath}
  \usepackage{bbm}
  \usepackage{cleveref}

   \crefformat{appendix}{#2}

  \newtheorem{theorem}{Theorem}[section]
  \newtheorem{lemma}[theorem]{Lemma}
  
  \newtheorem{proposition}[theorem]{Proposition}
  
  \newtheorem{definition}[theorem]{Definition}
  \newtheorem{corollary}[theorem]{Corollary}
  \newtheorem{remark}[theorem]{Remark}

  \def \E {{\mathbb E}}

  \newcommand{\cA}{{\cal A}}
  \newcommand{\cB}{{\cal B}}
  
  \newcommand{\cC}{{\cal C}}

  \newcommand{\cK}{{\cal K}}

  \newcommand{\cL}{{\cal L}}

  \newcommand{\cD}{{\cal D}}
  \newcommand{\cE}{{\cal E}}
  
  \newcommand{\cP}{{\cal P}}

  \def \ga {{\alpha}}
  \def \gb {{\beta}}
  \def \gep {{\varepsilon}}
  \def \gD {{\Delta}}
  
  \def \gl {{\lambda}}

  \def \gs {{\sigma}}
  
  \def \gga {{\gamma}}
  
  \def \go {{\omega}}
  \def \gd {{\delta}}

  \def \Z {{\mathbb Z}}
  \def \R {{\mathbb R}}
  
  \def \N {{\mathbb N}}
  
  \def \E {{\mathbf E}}
   \def \gP {{\mathbf P}}

  \def \gE {{\bf E}}

  \newcommand{\be}[1]{\begin{equation}\label{#1}}
  \newcommand{\ee}{\end{equation}}

  \newcommand{\bl}[1]{\begin{lemma}\label{#1}}
  \newcommand{\el}{\end{lemma}}

  \newcommand{\br}[1]{\begin{remark}\label{#1}}
  \newcommand{\er}{\end{remark}}

  \newcommand{\bt}[1]{\begin{theorem}\label{#1}}
  \newcommand{\et}{\end{theorem}}

  \newcommand{\bd}[1]{\begin{definition}\label{#1}}
  \newcommand{\ed}{\end{definition}}

  \newcommand{\bcl}[1]{\begin{claim}\label{#1}}
  \newcommand{\ecl}{\end{claim}}

  \newcommand{\bp}[1]{\begin{proposition}\label{#1}}
  \newcommand{\ep}{\end{proposition}}

  \newcommand{\bc}[1]{\begin{corollary}\label{#1}}
  \newcommand{\ec}{\end{corollary}}

  \newcommand{\bi}{\begin{itemize}}
  \newcommand{\ei}{\end{itemize}}

  \newcommand{\ben}{\begin{enumerate}}
  \newcommand{\een}{\end{enumerate}}

  \def \qed  {{\square\hfill}}
  \def \Z {{\mathbb Z}}
  \def \R {{\mathbb R}}
  
  \def \N {{\mathbb N}}
  
  \def \bE {{\mathbf E }}
  \def \ind {{\mathbf 1}}

\begin {document}

  \author{Julien Sohier\corref{cor1}\fnref{fn1}}
\ead{jusohier@gmail.com}
\address{Technische Universiteit Eindhoven, P.O. Box 513, 5600 MB Eindhoven, The Netherlands.}

      \begin{frontmatter}
 
 \title{ The scaling limits of the non critical strip wetting model.}
    
  \begin{abstract}
    The strip wetting model is defined by giving a (continuous space) one dimensional random walk $S$ a reward $\gb$ each
time it hits the strip $\R^{+} \times [0,a]$ (where $a$ is a positive parameter), which plays the role of a defect line. We show that this model exhibits a phase transition between a delocalized
regime ($\gb < \gb_{c}^{a}$) and a localized one ($\gb > \gb_{c}^{a}$), where the critical point $\gb_{c}^{a} > 0$ depends on $S$ and on $a$.
In this paper we give a precise pathwise description of the transition,
extracting the full scaling limits of the model. Our approach is based on Markov renewal theory. 
    \end{abstract}
    

  
\end{frontmatter}

  \section{Introduction and main results}
    
  \subsection{Definition of the models}
  
    We consider $ (S_n)_{n \geq 0}$ a random walk such that $S_0 := 0$ and $S_n := \sum_{i = 1}^n X_i$ where 
  the $X_i$'s are i.i.d. and $X_1$ has a density $h(\cdot)$ with respect to the Lebesgue measure. 
    We denote by $\gP$ the law of $ S$, and by $\gP_x$ the law of the same process starting from $x$.
    We assume that $h(\cdot)$ is continuous and bounded on $\R$, that $h(\cdot)$ is positive in a neighborhood of the origin,
    that $\gE[X] = 0$ and that $\gE [X^2]=: \gs^2 \in (0,\infty)$. We fix $a > 0$  in the sequel.
  
      The fact that $h$ is continuous and positive in the neighborhood of the origin entails that 
  \begin{equation}\label{hypn}
    n_0 := \inf_{n \in \mathbb{Z}^{+}} \left\{ (\gP[S_{n} > a],\gP[-S_{n} > a]) \in (0,1)^{2} \right\} < \infty.
  \end{equation} 
   
%
  
    For $N$ a positive integer, we consider the event $\cC_N := \{ S_ 1 \geq 0, \ldots, S_N \geq 0 \}$.
  We define the probability law (the \textit{free wetting model in a strip}) $\gP^f_{N,a,\gb}$ on $\R^N$ by 
  \begin{equation}
    \frac{d\gP^f_{N,a,\gb}}{d\gP} := \frac{1}{Z^f_{N,a,\gb}} \exp\left( \gb \sum_{k=1}^N \ind_{ S_k \in [0,a] } \right) \ind_{ \cC_{N}} 
  \end{equation}
  where $N \in \N$, $\gb \in \R$ and $Z^f_{N,a,\gb}$ is the normalization constant usually called the partition function of the system. 
  The second model we define is the \textit{constrained counterpart} of the above, that is 
  \begin{equation}
  \frac{d\gP^c_{N,a,\gb}}{d\gP} := \frac{1}{Z^c_{N,a,\gb}} \exp\left( \gb \sum_{k=1}^N \ind_{ S_k \in [0,a] } \right)
  \ind_{ \cC_{N}} \ind_{S_N \in [0,a]}. 
  \end{equation}
   
  Note in particular that 
  \begin{equation}
    \gP^c_{N,a,\gb} = \gP^f_{N,a,\gb}\left[ { } \cdot { } | S_N \in [0,a] \right],
  \end{equation}
      that $\gP^f_{N,a,0}$ is the law of $(S_1,\ldots, S_N)$ under the constraint 
  $\mathcal{C}_N := \{ S_1 \geq 0, \ldots, S_N \geq 0 \}$, 
  and that $\gP^c_{N,a,0}$ is the law of the same vector under the additional constraint $S_N \in [0,a]$.

    The process defined by the law $\gP_{N,a,\gb}$ is a  $(1+1)-$dimensional model for a linear chain of length $N$ which is 
  attracted or repelled to a defect \textit{strip} $[0,\infty) \times [0,a]$. By $(1+1)-$dimensional, we mean that the configurations of the linear
  chain are described by the trajectories $(i,S_{i})_{i \leq N}$ of the walk, so that we are dealing
  with directed models. The strength of this interaction with the strip is tuned by 
  the parameter $\gb$. Regarding the terminology, note that the use of the term \textit{wetting} has become customary to describe 
  the positivity constraint $\cC_N$ and refers to the interpretation of the field as an effective model for the interface of separation 
  between a liquid above a wall and a gas, see \cite{DGZ}. 

  It is an interesting problem to understand when the reward $\gb$ is strong enough to pin the chain near the defect strip, a 
  phenomenon that we call \textit{localization}, and what are the macroscopic effects of the reward on the system. In 
  this paper, we choose to characterize these effects through the scaling limits 
    of the laws $ \gP^c_{N,a,\gb}$ and $ \gP^f_{N,a,\gb}$. 
  More precisely, we first show the
    existence of a critical point $\gb_{c}^{a} > 0$ depending on $a$. This critical point separates 
     two phases: the localized phase and the delocalized one. Then we solve
    the \textit{full scaling limits} of the system in the case where $\gb \neq \gb_{c}^{a}$.

    We point out that these questions have been answered 
  in depth in the case of the
  standard wetting model, that is formally in the $a = 0$ case (see in particular \cite{IY} and \cite{DGZ}), and that 
  extending these results to our setup was 
    an open problem which has been raised by Giacomin (\cite[end of Chapter 2]{GB}). We stress that the
    techniques used in the standard wetting 
     model inspired most of the techniques used here, but that we had to overcome a number of technical problems.

  \subsection{The free energy.}
  
   A standard way to define localization for our models is by looking at the Laplace 
  asymptotic behavior of the partition function $Z^c_{N,a,\gb}$ as $N \to \infty$. More precisely,
  we define the free energy $F^{a}(\gb)$ by
  \begin{equation}\label{DFR}
    F^{a}(\gb) := \lim_{N \to \infty} \frac{1}{N} \log\left(Z^f_{N,a,\gb} \right)
  \end{equation}
  where the existence of the limit follows from Theorem \ref{delest} and Proposition \ref{esti}; we stress that the proof of these results 
   crucially relies on the representation
   which is made explicit in  Proposition \ref{ExFree}. 

    One basic observation is the fact that the free 
  energy is non-negative. 
  The following inequality holds:
    \begin{equation}
    \begin{split}
    Z^f_{N,a,\gb} & \geq  \bE\left[ \exp\left(\gb \sum_{k=1}^N \ind_{ S_k \in [0,a] } \right) \ind_{ S_k > a, k=n_{0} \ldots N}\right] \\
	& \geq \gP\left[ S_j > a, j = n_{0} \ldots, N\right]. \\
    \end{split}
    \end{equation} 
    
    For notational convenience, through the whole paper, for any event $B$ which is measurable with respect to $(S_{1},\ldots,S_{k})$, we will use
 the following notation:
   \begin{equation}
     \frac{1}{dx} \gP[ B, S_k \in dx] =: \gP[ B, S_k = x].
   \end{equation} 
    
      Then, integrating over $S_{n_{0}}$, we get:
  \begin{equation}
    \gP[ S_j > a, j = n_{0} \ldots, N] \geq \int_{(a,\infty)} \gP\left[S_{n_{0}} = t \right] \gP_t \left[ S_{1} > a, \ldots, S_{N-n_{0}} > a \right] dt.
  \end{equation} 
    
    We prove then in Lemma \ref{stayabove} below that 
    for some fixed $M$, the quantity $ N^{1/2}\gP_t \left[ S_1 > a, \ldots, S_{N-n_{0}} > a \right]  \in [c,c']$ for every $N \in \N$ and every
  $t \in [a,M]$, where $c,c'$ are positive constants. Thus:
    \begin{equation}
    Z^f_{N,a,\gb} \geq  \frac{c}{N^{1/2}}\int_{[a,M]} \gP\left[S_{n_{0}} = t \right] dt. 
    \end{equation} 
     
    Therefore $F^{a}(\gb) \geq 0$ for every $\gb$. Since the lower bound has been obtained by ignoring the contribution of the paths 
  that touch the strip, one is led to the following:

  \begin{definition}
  For $g \in \{ c,f\}$, the model $\{ \gP^g_{N,a,\gb} \}$ is said to be localized if $F^{a}(\gb) > 0$. It is said to be delocalized otherwise. 
  \end{definition}

    It is standard that $F^{a}(\cdot)$ is a convex increasing function, and in particular it is a continuous function as long as it is finite. 
    Therefore, there exists 
  a critical value $\gb_c^a \in \R$ such that the strip  wetting model is localized for $\gb > \gb_c^a$.  

  We stress that the terminology will be self-explanatory considering the scaling limits of these laws for 
   different $\gb'$s (see Theorem \ref{MAIN}). A more direct insight is to consider the 
   quantity
    \begin{equation}
     \frac{1}{N} \frac{\partial}{\partial \gb}\log(Z^f_{N,a,\gb}) = \gE^f_{N,a,\gb}\left[ \frac{1}{N} \sum_{j=0}^{N} \ind_{S_{j} \in [0,a]}\right],
    \end{equation}
    which, from standard convexity considerations, converges almost everywhere
    towards a positive quantity as soon as the model is localized, and vanishes in the limit $N \to \infty$ for $\gb < \gb_{c}^{a}$.

  \subsection{Scaling limits of the model.  }  
  
      We define the map $X^N: \R^N \mapsto C([0,1])$ (where  $C([0,1])$ is the 
  space of real continuous functions defined on $[0,1]$):
  \begin{equation} \label{X}
    X_t^N(x) := \frac{x_{\lfloor Nt \rfloor}}{ \gs N^{1/2}} + (Nt - \lfloor Nt\rfloor)
  \frac{x_{\lfloor Nt \rfloor+1} - x_{\lfloor Nt \rfloor}}{ \gs N^{1/2}}; t \in [0,1]
  \end{equation}
  where $\lfloor Nt \rfloor$ denotes the integer part of $Nt$. 
  Note that for a vector $x \in \R^{N}$, $X_t^N(x)$ is the linear interpolation of the process
  $\{ x_{ \lfloor Nt \rfloor}/ \gs N^{1/2} \}_{t \in \N/N \cap [0,1]}$. When the vector $x$ is random and distributed according to 
   the vector $S = (S_{1},\ldots,S_{N})$, this defines
    \begin{equation}\label{def3S}
     S^{N} = X^{N}(S)
    \end{equation} 
    an element of $C([0,1])$. Then we define the measures
  \begin{equation}
    Q_{N,a,\gb}^c := \gP_{N,a,\gb}^c \circ (X^N)^{-1} 
  \end{equation}
    and in an analogous way $Q_{N,a,\gb}^f$. These measures are defined on $C([0,1])$. In words, $Q_{N,a,\gb}^c$ (respectively $Q_{N,a,\gb}^f$) is the law of 
   $S^{N}$ when $S$ is distributed according to $\gP_{N,a,\gb}^c$ (respectively $\gP_{N,a,\gb}^f$). 
  
    We consider the following standard processes:
  \begin{enumerate}
    \item[$\star$]  the Brownian motion $(B_t)_{t \in [0,1]}$.
    \item[$\star$]  the Brownian meander  $(m_t)_{t \in [0,1]}$ which is the Brownian motion conditioned to stay positive on $[0,1]$.
    \item[$\star$]  the normalized Brownian excursion $(e_t)_{t \in [0,1]}$ which is the brownian bridge
  conditioned to stay positive on $[0,1]$. 
  \end{enumerate}

    Our main result is the following:
  \begin{theorem} \label{MAIN}
    Both the free and the constrained models undergo a wetting transition at $\gb = \gb^a_c$. 
  More precisely: 
  \begin{enumerate}
    \item  in the subcritical regime, that is if $\gb < \gb^a_c$, then
    \begin{itemize} 
    \item $(Q_{N,a,\gb}^c)_N$ converges weakly in $C([0,1])$ to the law of $e$.
    \item $(Q_{N,a,\gb}^f)_N$ converges weakly in $C([0,1])$ to the law of $m$.
    \end{itemize}
    \item in the supercritical regime, that is if $\gb > \gb^a_c$, then both $(Q_{N,a,\gb}^c)_N$ and $(Q_{N,a,\gb}^f)_N$ converge 
    in $C([0,1])$ to the measure concentrated on the constant function taking value zero.
    \end{enumerate}
  \end{theorem}
  
   The corresponding result in the case of the standard wetting model has been shown in \cite[Theorem 1]{DGZ}. According to the results 
    proved in their setup, we expect that at the critical point, the limiting process should be $(|B_{t}|)_{t \in [0,1]}$ in the free case 
     and $e$ in the constrained case. We stress 
      that these results have been shown in a weak sense (that is only at the level of contact sets) in \cite[Theorem 1.5]{Soh1}, and that extending 
       them to the convergence in law at the level of processes would require some additional work, in particular to deal 
        with tightness issues in this case. We also expect Theorem \ref{MAIN} to hold also in the case of a random walk in the domain of attraction 
         of an $\ga$ stable law, at the cost of several technicalities though; for example, one should deal with the lack of regularity of the limiting 
         processes involved (see Section \ref{TDP} for details). 
        
        Given 
   a vector $(x_{1},\ldots,x_{N}) \in \R^{N}$, we define
   \begin{equation}
    \cA(x) :=  \left\{ i, x_{i} \in [0,a] \right\}.
   \end{equation} 
  
  The following result is a crucial step for the proof of the first part of Theorem \ref{MAIN} and is interesting in itself. 
   It states that in the subcritical phase, the "dry region" (that is the set $\cA$)
   reduces to a finite number of points all 
  being at a finite distance from $\{ 0 \}$ in the free case, from $\{ 0 \}$ and from $\{N\}$ in the constrained case. It
  is the analogous of \cite[Proposition 5]{DGZ}. 
  
  \begin{theorem} \label{TRLO}
    For $\gb < \gb^a_c$, the following convergences hold:
    \begin{equation}
			  \lim_{ L \to \infty} \limsup_{ N \to \infty} \gP_{N,\gb,a}^f\left[ \max \cA \geq L \right] = 0 \\
    \end{equation}
  and
  \begin{equation} \label{CCCCC}
    \begin{split}
		      & \lim_{ L \to \infty} \limsup_{ N \to \infty} \gP_{N,\gb,a}^c\left[ \max ( \cA \cap [1, N/2]) \geq L \right] = 0,\\
		    & \lim_{ L \to \infty} \limsup_{ N \to \infty} \gP_{N,\gb,a}^c\left[ \min ( \cA \cap [N/2,N]) \leq N-L \right] = 0. \\
		    \end{split}
  \end{equation}
  \end{theorem}
   

  Theorem \ref{MAIN} characterizes the Brownian scaling of the model when $\gb \neq \gb^a_c$. Infinite scaling results like Theorem \ref{MAIN} 
  have been proved in different contexts involving polymer measures. The first mathematical paper dealing with such an issue is \cite{IY}
  where the authors proved an analogous convergence in the homogeneous pinning model 
  for the case where $S$ is a symmetric random walk with increments taking values in $\{-1,0,1 \}$. Their results have
  been strongly generalized in \cite{DGZ} where the same assumptions are made on $S$ as in this paper, and a further 
  generalization of their results in the case where $S$ is in the domain of attraction of the standard normal law has been obtained in
  \cite{CGZ}.
  
    Analogous results have also been obtained in \cite{CGZ1} in the case of inhomogeneous, but periodic pinning models, and 
    more recently in \cite{CD2} in the case where the interaction is 
  of \textit{Laplacian} type. Related models with 
      a different characterization of the large scale limits have been
      considered recently \cite{FuN}. Finally,  a closely related pinning model in continuous time has been considered 
  and resolved 
  in \cite{CrKMV}; we stress however that their techniques are very peculiar to the continuous time setup.

  \subsection{ Organization of the paper } 

  To prove our main result Theorem \ref{MAIN} in the localized phase, the main point has been to show 
   a general state space Markov renewal theorem (which is Theorem \ref{MarRen}).
   
     We stress that dealing with the delocalized phase is technically much heavier. The first observation is the fact that
   a common feature shared by the strip wetting model and the classical homogeneous one is the fact that the measures $ \gP^{g}_{N,a,\gb}$
   (for $g \in \{c,f \}$) exhibit
  a remarkable decoupling between the contact level set $\mathcal{A} = \{ i \leq N, S_i \in [0,a] \}$ and the excursions of $S$ between 
  two consecutive contact points. More precisely, conditioning on 
    $\cA = \{ t_1, \ldots, t_k\}$ and on $(S_{t_1}, \ldots, S_{t_k})$, the \textit{bulk} excursions 
  $e_i = \{ e_i(n) \}_n := \left\{ \{S_{t_i+n} \}_{ 0 \leq n \leq t_{i+1} - t_i} \right\}$
  are independent under  $ \gP^{g}_{N,a,\gb}$ and are distributed like the walk $(S',\gP_{ S_{t_i}})$ conditioned on the 
  event
   \begin{equation}
    \left\{ S'_{t_{i+1} - t_i}  = S_{t_{i+1}} , S'_{t_i + j} > a, j \in \{ 1, \ldots, t_{i+1} - t_i -1 \} \right\}.
   \end{equation}
    
  It is therefore clear that 
  to extract the scaling limits of the laws $\gP^{g}_{N,a,\gb}$, one has to combine good control over the law of the contact set $\mathcal{A}$ 
  and suitable asymptotics properties of the excursions. This decoupling turns out to be the starting point of our proofs, see Section \ref{TDP} for details. 
  
   More precisely, here is the plan of this paper: 
  \begin{enumerate}
    \item[-] in Section \ref{ShFT}, we first give some 
   recent local limit estimates for random walks conditioned to stay non negative (see Lemma \ref{fluctu}), and we use them to prove a local 
    limit estimate related to our problem (Lemma \ref{stayabove}); then we 
    recall a result on the tails of the return probability  
  to the strip for large $N$ (Lemma \ref{Pr}), which has been proved in \cite[Theorem 3.1]{Soh1}. 
    \item[-] in Section \ref{MRT}, we give a representation of $F^{a}(\cdot)$ and of $\gb_{c}^{a}$ in terms of the spectral radius 
    of a Hilbert Schmidt operator (Proposition \ref{ExFree}); then 
     we show that the set of contact points with the strip under $\gP_{N,a,\gb}^c$ is distributed according to the law 
  of a Markov renewal process conditioned to hit the strip at time $N$ (Proposition \ref{RMRP}). This representation implies a very useful 
  expression for the partition function $Z_{N,a,\gb}^c$ which is the key to our main results. 
    \item[-] in Section \ref{TLP}, we deal with the localized phase and we 
    make use of a finite mean Markov renewal theorem (Theorem \ref{MarRen}) to deduce asymptotic estimates on $Z_{N,a,\gb}^c$ and $Z_{N,a,\gb}^f$ 
    (which are given in Theorem \ref{delest}). These estimates are enough to prove Theorem \ref{MAIN} in the localized phase.
  \item[-]  Section \ref{TDP} is devoted to the proof of Theorems \ref{TRLO} and \ref{MAIN} in the localized phase. These proofs are carried out 
  by first giving estimates on the partition functions (see Proposition \ref{esti}) by the means of an (infinite mean) Markov renewal theorem (which 
   has been proved in \cite[Section 7.2]{CD1}); 
  from these estimates we deduce Theorem \ref{TRLO}. Since the process conditioned on the contact set behaves like the free random walk, we 
    can then combine these results with powerful limit theorems which have been obtained 
  in \cite{Shi} (for the free case) and much more recently in \cite{CaCha} (for the constrained case) to obtain Theorem \ref{MAIN}.   
  \item[-] in Appendix \ref{apA}, we show the Markov renewal theorem (Theorem \ref{MarRen}) in the finite mean case in a general 
  framework. In  Appendix
  \ref{valcr}, we illustrate the power of Proposition \ref{ExFree}
  in the particular case of the $(p,q)$ random walk for $a=1,2$ (which is Theorem \ref{caspq}).
  \end{enumerate}

  \section{Preliminary facts}\label{ShFT}  
  
  \subsection{Recurrent notations and terminology}  
  
  For $a,b \in \R$, we define $a \vee b := \max(a,b)$.
   
    For $a_n,b_n$ positive sequences, we write $a_n \sim b_n$ if $\lim_{n \to \infty} a_n/b_n = 1$. By a slight abuse 
     of notation, we also write  $a_n \sim b_n$ in the case where 
    the sequences $a_{n}$ and $b_{n}$ are identically null. 
     
    More generally, for $a_n(x)$ a positive sequence depending on a
  parameter $x\in \gD$ where $\gD$ is a subset
  of $\R^d, d \geq 1$, $\ga \in \R$ and $b(\cdot)$ 
  a measurable function on $\gD$,  we often say that the equivalence 
  \begin{equation}
    a_n(x) \sim \frac{b(x)}{n^{\ga}}
  \end{equation}
  holds \textit{uniformly for $x$ in $\gD$} if the following holds:
  \begin{equation}
      \lim_{n \to \infty} \sup_{x \in \gD} | n^{\ga} a_n(x) - b(x)| = 0.
  \end{equation}

  In this paper, we deal with kernels of two kind. Kernels of the first kind are just $\gs$-finite kernels on $\R$, that 
  is functions $A: \R \times \cB(\R) \mapsto \R^+$ (where $\cB(\R)$ denotes the Borel $\gs$-field of $\R$), 
  and such that for each $x \in \R$, 
  $A_{x,\cdot}$ is a $\gs$-finite measure on $\R$ and $A_{\cdot,F}$ is a Borel function for every $F \in \cB(\R)$.
  Given two such kernels 
  $A$ and $B$, their composition is denoted by $(A \circ B)_{x,dy} := \int_{z \in \R} A_{x,dz}B_{z,dy}$ and of course 
  $A^{\circ k}_{x,dy}$ denotes the $k$-fold composition of $A$ with itself where  $A^{\circ 0}_{x,dy} := \gd_x(dy)$. 

    The second kind of kernels is obtained by letting a kernel of the first kind depend on a further parameter $n \in \Z^+$: more precisely, 
  we consider objects of the form $A_{x,dy}(n)$ with $x,y \in \R, n \in \Z^+$. Given two such
  kernels $A_{x,dy}(n), B_{x,dy}(n)$ we define
  their convolution 
  \begin{equation}
    (A \ast B)_{x,dy}(n) := \sum_{m=0}^{n} (A(m) \circ B (n-m))_{x,dy} = \sum_{m=0}^{n} \int_{\R} A_{x,dz}(m) B_{z,dy} (n-m),
  \end{equation}
    and the $k$-fold convolution of the kernel $A$ with itself is denoted by $A^{\ast k}_{x,dy}$ where by definition 
    $A^{\ast 0}_{x,dy} := \gd_x(dy) \ind_{ n = 0}$. Finally given two kernels $A_{x,dy}(n)$ and $B_{x,dy}$ and a positive sequence $a_n$, 
  we write 
  \begin{equation}
    A_{x,dy}(n) \sim \frac{B_{x,dy}}{a_n} 
  \end{equation}
    to mean $ A_{x,F}(n) \sim \frac{B_{x,F}}{a_n} $ for every $x \in \R$ and for every bounded set $ F \subset \R$.   
    
      Natural and useful examples of the above kernels are the partition functions; 
      namely, for $x,y \in [0,a]\times \R^{+}$, we define:
  \begin{equation}
  Z^c_{N,a,\gb}(x,dy) := \bE_x\left[ \exp\left( \gb \sum_{k=1}^N \ind_{ S_k \in [0,a] } \right)
  \ind_{ S_k \geq 0, k=1 \ldots N} \ind_{S_N \in dy} \ind_{ y \in [0,a]} \right],
  \end{equation}
  and its free counterpart 
    \begin{equation}
  Z^f_{N,a,\gb}(x,dy) := \bE_x\left[ \exp\left( \gb \sum_{k=1}^N \ind_{ S_k \in [0,a] } \right)
  \ind_{ S_k \geq 0, k=1 \ldots N} \ind_{S_N \in dy} \right].
  \end{equation}

    \subsection{ Markov renewal and random walk fluctuation theory.} \label{MRSS}

    Let us introduce the following transition kernel:  
    \begin{equation} 
  \begin{split}
		& F_{x,dy}(n) := \gP_{x} [S_1 > a, S_2 > a, \ldots, S_{n-1} > a, S_n \in dy] \ind_{ x,y \in [0,a]} \hspace{.2 cm} \text{if} \hspace{.2 cm}  n \geq 2, \\
		&  F_{x,dy}(1) := h(y-x) \ind_{ x,y \in [0,a]} dy.
		    \end{split}
    \end{equation}
    
    We write $f_{x,y}(n)$ for the density of $F_{x,dy}(n)$ with respect to the Lebesgue measure.

    We denote by $(\tau_n)_{n \geq 0}$ the times of return to $[0,a]$ of $S$, that is $\tau_0 := 0$ and,
  for $n \geq 1$, $\tau_{n} := \inf \{ k> \tau_{n-1} | S_k \in [0,a] \}$.
  Note that $(\tau_n)_{n \geq 0}$ is not a true renewal process. Introducing the process $(J_n)_{n \geq 0}$ where  $J_n := S_{\tau_n}$, the
  process $\tau$ is a so called \textit{Markov
  renewal process} whose modulating chain is the Markov chain $J$. The topic of Markov renewal theory is a classical one, a well known 
  reference is \cite{As}. 
    
  We finally  denote by $l_N$ the cardinality of $\{ k \leq N | S_k \in [0,a] \}$. 
  With these notations, we can write the joint law of $( l_N, (\tau_n)_{  n \leq l_N}, (J_n)_{n \leq l_N})$
  under $\gP^c_{N,a,\gb}$ under the following form: 
    \begin{equation}\label{HPP}
  \begin{split}
    &  \gP^c_{N,a,\gb} [ l_N = k,\tau_i = t_i, J_i \in dy_i, i=1 ,\ldots, k ] \\
  & \phantom{xxxxxx}  = \frac{e^{\gb k}}{Z^c_{N,a,\gb}}  F_{0,dy_1}(t_1)F_{y_1,dy_2}(t_2-t_1) \ldots F_{y_k,dy_{k-1}}(N-t_{k-1}) 
  \end{split}
      \end{equation} 
  where $k \in \N,  0 < t_1 < \ldots < t_k =N$ and $ (y_i)_{i=1, \ldots,k} \in \R^k$.

    It is then clear that getting asymptotic estimates on the partition functions $Z^c_{N,a,\gb}$ (and thus $Z^f_{N,a,\gb}$),
  requires an accurate 
  control on the asymptotic behavior of $F_{\cdot,\cdot}(n)$ for large $n$.

  To achieve this, we collect some basic facts about random walk fluctuation theory.

  For $n$ an integer, we denote by $T_n$ the {\sl $n^{\text{th}}$ ladder epoch}; that is $T_0 := 0$
  and, for $n \geq 1$, $T_n := \inf\{ k \geq T_{n-1}, S_k > S_{T_{n-1}} \}$. 
  We also introduce the so-called {\sl ascending ladder heights} $ (H_n)_{n \geq 0}$,
  which, for $ k \geq 1$, are given by $H_k := S_{T_k}$. Note that the process $(T,H)$ is 
  a bivariate renewal process on $(\R^+)^2$. 
    In a similar way, we write $(T^-,H^-)$ for the strict descending ladder variables process, which is defined by
    $(T^-_0,H^-_0) := (0,0)$ and 
    \begin{equation}
      T^-_n := \inf\{ k \geq T_{n-1}, S_k < S_{T_{n-1}} \} \hspace{.6 cm} \text{and} \hspace{.6 cm} H^-_k := -S_{t_k^-}. 
    \end{equation}
    
     Let us consider the renewal function $U(\cdot)$
 associated to the ascending ladder heights process:
 \begin{equation}
 U(x) := \sum_{k=0}^{\infty} \gP [ H_k \leq x] = \bE[ \mathcal{N}_x] = \int_{0}^{x} \sum_{m=0}^{\infty} u(m,y) dy
\end{equation}
  where $\mathcal{N}_x$ is the cardinality of $\{ k \geq 0, H_k \leq x \}$ and
 $ u(m,y) := \frac{1}{dy} \gP[ \exists k \geq 0, T_{k} = m, H_{k} \in dy]$ is the renewal mass function associated to $(T,H)$.
 It follows in particular from this definition 
 that $U(\cdot)$ is a subadditive increasing function, and in our context it is also continuous. Note also that $U(0) = 1$. 
  We denote by $V(x)$
  the analogous quantity for the process $H^-$, and by $v(m,y)$ the renewal mass function associated
 to the descending renewal $(T^{-},H^{-)}$. 
 
   The following local limit estimates have been proved recently (\cite{Do3} and \cite{CaCha}):
   \begin{lemma}\label{fluctu} Uniformly on sequences $x_{n}, y_{n}$ such that $x_{n} \vee y_{n} = o(\sqrt{n})$, the following 
    equivalences hold:
     \begin{equation}\label{pacons}
      \gP_{x_{n}}[S_{1} \geq 0, \ldots, S_{n} \geq 0 ] \sim V(x_{n}) \gP[T_{1}^{-} > n] \sim \frac{V(x_{n}) }{ \sqrt{2 \pi } \gs \sqrt{n}}
     \end{equation} 
      and
      \begin{equation}
       \gP_{x_{n}}[ S_{1} \geq 0, \ldots, S_{n} \geq 0, S_{n} = y_{n}] \sim  \frac{V(x_{n})U(y_{n})}{n} \gP[S_{n} = y_{n}].  
      \end{equation} 
   \end{lemma}

      Note that making use of Gnedenko's classical local limit theorem, for sequences $(x_{n}),(y_{n})$ satisfying
      the same assumptions as in Lemma \ref{fluctu}, one gets the equivalence
      \begin{equation}\label{const}
               \gP_{x_{n}}[ S_{1} \geq 0, \ldots, S_{n} \geq 0, S_{n} = y_{n}] \sim \frac{V(x_{n})U(y_{n})}{\gs \sqrt{2 \pi}n^{3/2}}.  
      \end{equation}
       
       The following result is a consequence of Lemma \ref{fluctu}:
       \begin{lemma}\label{stayabove}
         For any $x \in [0,a]$, one has the following convergence: 
        \begin{equation}
         \gP_{x}[ S_{1} > a, \ldots, S_{n} > a] \sim \frac{\gP[H_{1} \geq a - x]}{\sqrt{2 \pi } \gs n^{1/2}}. 
        \end{equation}   
        
       \end{lemma}
        
        Note that both terms in the above equivalence might be identically $0$ (at least for $x=0$ in the case $n_{0} > 1$), so that 
         we recall that we use the convention $0 \sim 0$ 
          and we note the equivalence (valid 
         for any $u > 0$)
       \begin{equation}
        \gP[H_{1} > u] > 0 \Longleftrightarrow P[S_{1} > u] > 0. 
       \end{equation} 
       
    \textit{Proof of Lemma \ref{stayabove}} 
     We integrate over $S_{1}$ to get:
     \begin{equation}
     \begin{aligned}
       \gP_{x}[ S_{1} > a, \ldots, S_{n} > a]  =  \int_{ u \in [a, \infty)} \gP_{x}[ S_{1} = u, \ldots, S_{n} > a] du \\
                  = \int_{ u \in [a, \infty)} h(u-x) \gP_{u-a}[ S_{1} > 0, \ldots, S_{n-1} > 0] du \\
       =  \int_{ u \in [a,  n^{1/4}]} h(u-x) \gP_{u-a}[ S_{1} > 0, \ldots, S_{n-1} > 0] du \\
     + \int_{ u \in [n^{1/4}, \infty)} h(u-x) \gP_{u-a}[ S_{1} > 0, \ldots, S_{n-1} > 0] du.
     \end{aligned}
     \end{equation} 
      
      For the second term in the right hand side of the above equalities, we immediately get that, for any $n$ large enough:
      \begin{equation}
       \begin{aligned}
    \int_{ u \in [n^{1/4}, \infty)} h(u-x) \gP_{u-a}[ S_{1} > 0, \ldots, S_{n-1} > 0] du & \leq  \int_{ u \in [n^{1/4}-a, \infty)} h(u) du \\ 
   & \leq  \frac{4}{ n^{1/2}}\int_{ u \in [n^{1/4}/2, \infty)} u^{2} h(u) du, 
       \end{aligned} 
          \end{equation} 
     and since $\gE[X^{2}] < \infty$, it immediately follows that this term is $o(n^{-1/2})$. 
      On the other hand, making use of Lemma \ref{fluctu}, we get that
       \begin{equation}
        \int_{ u \in [a,  n^{1/4}]} h(u-x) \gP_{u-a}[ S_{1} > 0, \ldots, S_{n-1} > 0] du \sim 
         \int_{ u \in [a,  n^{1/4}]} h(u-x) \frac{V(u-a) }{ \sqrt{2 \pi } \gs \sqrt{n}} du. 
       \end{equation} 
       
        Then we recall that, using duality arguments (see for example \cite[Proof of Theorem 3.1]{Soh1} for a proof), one can show that
      \begin{equation}
       \int_{ u \in [a, \infty)} h(u-x) V(u-a) du = \gP[H_{1} \geq a - x],
      \end{equation} 
       from which we finally deduce Lemma \ref{stayabove}. 
     $\qed$

  We define the following function: 
    \begin{equation}\label{DefPh}
    \Theta_a(x,y) :=  \frac{\gP [H_1^- \geq a-y]\gP [H_1 \geq a-x]}{\gs \sqrt{2\pi} } \hspace{ .2 cm}
    \ind_{ x,y \in [0,a]}.
    \end{equation}
    
  By using similar techniques as the ones we just developed for the proof 
   of Lemma \ref{stayabove}, in \cite{Soh1} we showed the following result, which is the cornerstone of our approach: 
    \begin{lemma} \label{Pr} The following equivalence holds uniformly on $(x,y) \in [0,a]^2$: 
      \begin{equation}\label{EQC}
      n^{3/2}f_{x,y}(n) \sim \Theta_a(x,y).
      \end{equation}
	\end{lemma}
	
    A similar remark as after Lemma \ref{stayabove} holds here as well. 
	Since $\Theta_{a}$ is bounded on $[0,a]^{2}$, a trivial consequence of the above result is the fact that the left hand side in \eqref{EQC} is
	  dominated by a multiple of its right hand side.

  \section{ An infinite dimensional problem} \label{MRT}

    \subsection{Defining the free energy} In this Section, we define the free energy in a
    way that allows us to make use of the Markov renewal structure that we pointed at in the previous Section. 
    For $\gl  \geq 0$, we introduce the following kernel:
  \begin{equation}\label{defker}
  B_{x,dy}^{\gl} := \sum_{n=1}^{\infty} e^{-\gl n} F_{x,dy}(n)
  \end{equation}
  and the associated integral operator 
  \begin{equation}
    (B^{\gl}h)(x) := \int_{[0,a]} B_{x,dy}^{\gl} h(y). 
  \end{equation} 
  
   Making use of the asymptotics \eqref{EQC}, one can show, as in \cite[Lemma 4.1]{CD1}, that for any $\gl \geq 0, B_{x,dy}^{\gl}$ is a compact operator on
   the Hilbert space $L^2([0,a])$. Using this fact, we introduce $\gd^a(\gl)$, the spectral radius of the operator $B^{\gl}$, which is an isolated
    and simple eigenvalue of  $B_{x,dy}^{\gl}$ (see Theorem 1 in \cite{Ze}). The function $\gd^a(\cdot)$ is non-increasing,
  continuous on $[0,\infty)$ and analytic on $(0,\infty)$ because the operator $B_{x,dy}^{\gl}$ has these properties. The analyticity
  and the fact that $\gd^a(\cdot)$ is not constant (as $\gd^a(\gl) \to 0$ as $\gl \to \infty$) force $\gd^a(\cdot)$ 
  to be strictly decreasing.
    
    We denote by $(\gd^a)^{-1}(\cdot)$ its
    inverse function, defined on $(0,\gd^a(0)]$. We are now ready to state the following fundamental Proposition.

   \begin{proposition}\label{ExFree}
    We have the equalities:
      \begin{equation} \label{DefFEE}
    \gb_c^a := -\log(\gd^a(0)), \hspace{.2 cm} F^a(\gb) := ( \gd^a)^{-1}(\exp(-\gb) ) \hspace{.2cm} \text{if} \hspace{.2cm} \gb \geq \gb_c^a  \hspace{.2cm} \text{and} \hspace{.2cm}  0 \hspace{.2cm} \text{otherwise.}
  \end{equation}
   \end{proposition}
   
    The equalities in Proposition \ref{ExFree} are direct consequences of Theorem \ref{delest} and Proposition \ref{esti} which are proved in the following parts in a
     way which does not depend on Proposition \ref{ExFree} (as will be clear from their proofs). We stress that they are central tools 
      in the standard wetting model, in particular to compute exactly the critical point and to show the critical behavior of the free energy. 
      We refer to 
       \cite[Chapter 2]{GB} for more details. We illustrate the power of this Proposition 
        in Appendix \ref{valcr} for the particular case of the $(p,q)$ random walk.
  

    \subsection{A useful representation for  $Z^c_{N,a,\gb}$ }\label{UREP}
     
      For $x \in [0,a]$, we denote by $b(x,\cdot)$ the density of
      $B_{x,dy}^{F^{a}(\gb)}$ with respect to Lebesgue measure. Combining \eqref{defker}, Lemma \ref{Pr}
       and the positivity of $\Theta_{a}$ on $[0,a]^{2}$, 
       we deduce that $b(x,y) > 0$  for every $(x,y)\in [0,a]$. This fact implies the uniqueness (up to a multiplication 
  by a positive constant) and the positivity almost everywhere of 
  the right (respectively the left) Perron Frobenius eigenfunctions $v_{\gb}(\cdot)$ (respectively $w_{\gb}(\cdot)$) of $B_{x,dy}^{F^a(\gb)}$. We refer 
   to \cite[Section 4.2]{CD1} for more details. In particular, one can show that the function $v_{\gb}(\cdot)$ is positive everywhere (not only almost
    everywhere); hence we can define the kernel 
  \begin{equation}
    K_{x,dy}^{\gb}(n) := e^{\gb} F_{x,dy}(n) e^{-F^a(\gb)n}  \frac{v_{\gb}(y)}{v_{\gb}(x)},
  \end{equation}
     and using the definition \eqref{DefFEE} of $\gb_{c}^{a}$, it is easy to check that
    \begin{equation} \label{INVMP}
    \begin{split}
    & \int_{y \in \R} \sum_{n \in \N} K_{x,dy}^{\gb}(n) = \min\left(1, \frac{e^{\gb}}{e^{\gb_c^a}}\right). 
    \end{split}
    \end{equation}

  Then we define the law $\cP^{\gb}$ under which the joint
  process  $(\tau_k,J_k)_{k \geq 0}$ is an inhomogeneous Markov chain (defective if $\gb < \gb_c^a$) on $\Z^+ \times [0,a]$ by:
  \begin{equation} \label{Def}
    \cP^{\gb} \left[ (\tau_{k+1}, J_{k+1}) \in ( \{ n \} , dy)  | (\tau_{k}, J_{k}) = (m,x)  \right] := K^{\gb}_{x,dy} (n-m).
  \end{equation}
    
  The sequence $(\tau_k)_{k \geq 0}$ is a Markov renewal, the process $(J_i)_{i \geq 0}$ 
  being its modulating chain.
    We then have the following property, whose proof is contained in 
    \eqref{HPP}: 
  \begin{proposition} \label{RMRP}
    For any $N \in \N$, the vector $( l_N, (\tau_n)_{  n \leq l_N}, (J_n)_{n \leq l_N})$ has the same
  law under $\gP^c_{N,a,\gb}$ as under the conditional law $\cP^{\gb}(\cdot | N \in \tau)$. Equivalently:
  \begin{equation}
    \begin{split}
      & \gP^c_{N,a,\gb} \left[ l_N = k,\tau_j = t_j, J_j \in dy_i, i=1 ,\ldots, k \right]\\
    & \phantom{xxxxxx} = \cP^{\gb} \left[ l_N = k,\tau_j = t_j, J_j \in dy_i, i=1 ,\ldots, k | N \in \tau \right].
    \end{split}
  \end{equation}
  \end{proposition}
  
  Proposition \ref{RMRP} shows in particular that the partition  $Z^c_{N,a,\gb} $ can be interpreted as 
  the Green function associated to the Markov renewal $\tau$, that is
  $Z^c_{N,a,\gb} = \cP^{\gb}[ N \in \tau]$. More generally, 
  for $x,y \in [0,a], Z^c_{N,a,\gb}(x,dy) =  \cP^{\gb}[ \exists k, \tau_k = N, J_k \in dy | J_{0} = x]$. Equivalently, 
  we have the equality 
    \begin{equation} \label{Rep}
    Z^c_{N,a,\gb} = \exp(F^a(\gb)N) \int_{[0,a]} \frac{v_{\gb}(0)}{v_{\gb}(y)} \sum_{k \geq 0} (K^{\gb})^{\ast k}_{0,dy}(N)
  \end{equation}
  which is a consequence of the more general equality:
  \begin{equation} \label{FUN}
      Z^c_{N,a,\gb}(x,dy) = \exp(F^a(\gb)N)  \frac{v_{\gb}(x)}{v_{\gb}(y)} \sum_{k \geq 0} (K^{\gb})^{\ast k}_{x,dy}(N) 
  \end{equation}
  which holds for $x,y \in [0,a]$.

  \section{The localized phase} \label{TLP}

    Let $\gb > \gb_c^a$. 
      In this case, the two functions $w_{\gb}(\cdot)$ and $v_{\gb}(\cdot)$ are uniquely defined up to a multiplicative constant, 
  and we use this degree of freedom to fix $\int_{\R} v_{\gb}(x)w_{\gb}(x) \ind_{ x \in [0,a]} =1$. 
  Thus the measure $\mu_{\gb}$ defined by 
  \begin{equation}
    \mu_{\gb}(dx) := v_{\gb}(x)w_{\gb}(x) \ind_{ x \in [0,a]} dx
  \end{equation}
    is a probability measure. It is a straightforward computation to verify that for $\gb > \gb_c^a$, the probability  $\mu_{\gb}$ is 
  invariant for the kernel  $ \sum_{n \geq 1} K_{x,dy}^{\gb} (n) $, and hence for the Markov process $(J_n)$.

   The next result is a then consequence of Theorem \ref{MarRen}, which we show in Appendix \ref{apA}. 
    
    \begin{theorem} \label{delest} 
  For $\gb > \gb_c^a$,  for every $x \in [0,a], y \in \R^+$, as $N \to \infty$, one has the convergence  :
  \begin{equation}
    Z^c_{N,a,\gb}(x,dy) \sim \frac{ v_{\gb}(x) v_{\gb}(y)}{ C_{\gb}} \exp(F^a(\gb)N) dy
  \end{equation}
    where for a fixed $x \in [0,a]$,
  the convergence of $ Z^c_{N,a,\gb}(x,dy) \exp(-F^a(\gb)N)$ towards 
  $\frac{ v_{\gb}(x) v_{\gb}(y)}{ C_{\gb}}dy$ holds in total variation norm. 

    These estimates imply in particular that there exist two positive constants $C^a(\gb)$ and $C^a_f(\gb)$ such that, :
    \begin{enumerate}
  \item $ Z^c_{N,a,\gb} \sim C^a(\gb) \exp(F^a(\gb)N) $,
  \item $ Z^f_{N,a,\gb} \sim C^a_f(\gb) \exp(F^a(\gb)N)$.
  \end{enumerate}
  \end{theorem}
  
   This result readily implies a much finer result on the scaling limits of the system
    than the one of Theorem \ref{MAIN} in the localized phase, namely the fact that for $\gb > \gb_{c}^{a}$ the Markov Renewal process 
    defined in \eqref{Def} converges without need of rescaling, and in particular this last fact trivially implies Theorem \ref{MAIN} in the localized phase.
    We refer to \cite[Chapter 2]{GB} for more details.

  \textit{Proof of Theorem \ref{delest}}
  
   Of course, since for $\gb > \gb_c^a$,
   \begin{equation}
    C_{\gb} = \int_{(x,y) \in [0,a]^2} \mu_{\gb}(dx) \sum_{k \geq 1} k K^{\gb}_{x,dy}(k) < \infty,
   \end{equation} 
   a  direct consequence of (and our main motivation for proving) Theorem \ref{MarRen} is the fact that 
  in the localized regime, for any $x \in [0,a]$, the following convergence holds in total variation norm: 
    \begin{equation} \label{Eq111}
  \lim_{ N \rightarrow \infty} \cP^{\gb} \left[ \exists k \in \N, \tau_k = N, J_N \in dy  | J_{0} = x\right] = 
  \frac{ \mu_{\gb}(dy)}{ C_{\gb} }.
  \end{equation}
   
    Combining this with identity \eqref{FUN}, we get  :
  \begin{equation}
  Z^c_{N,a,\gb}(x,dy) \sim \exp(F^a(\gb)N)  \frac{v_{\gb}(x)}{v_{\gb}(y)} \frac{\mu_{\gb}(dy)}{ C_{\gb}}.
  \end{equation}
  
    The free case follows from the asymptotic behavior of $Z^c_{N,a,\gb}$ and the relation
    \begin{equation}
    Z^f_{N,a,\gb} = e^{F^a(\gb)N} \sum_{t=0}^N \int_{x \in [0,a]}  Z^{c}_{N-t,a,\gb}(dx) e^{-F^a(\gb)(N-t)} \gP_x[S_1 >a, \ldots, S_t >a] e^{-F^a(\gb)t}.
    \end{equation}
    
     Note that since $h(\cdot)$ is positive in the neighborhood of the origin, $\gP_x[S_1 >a, \ldots, S_t >a]$ is strictly positive at least for $x$ close enough 
      to $a$. 
    
    This entails:
    \begin{equation}
      Z^f_{N,a,\gb} \sim C_a(\gb) e^{F^a(\gb)N} \sum_{t=0}^{\infty} e^{-F^a(\gb)t} \int_{[0,a]} \frac{ \mu_{\gb}(dx)}{C_{\gb}} \gP_x[\tau_1 > t+1].
    \end{equation}

  $ \qed $

  \section{The delocalized phase} \label{TDP}
  
   \subsection{ Some results borrowed from the standard homogeneous wetting.}
  
   We stress that we can adapt in a straightforward way some of the techniques borrowed from different 
    papers on the topic of scaling limits linked to polymer models to our case of interest. We first mention that, using
    techniques which have been developped in \cite{CD1}, we can prove the
    following asymptotics on the partition functions in the delocalized phase: 
    
%
%
  \begin{proposition} \label{esti}
    For $\gb < \gb_c^a$,as $N \to \infty$, we have the following:
  \begin{enumerate}
  \item $ Z^c_{N,a,\gb}(x,dy) \sim C'^a(\gb) N^{-3/2} \Theta_{a}(x,y) dy $
  \item $ Z^f_{N,a,\gb}(x,dy) \sim C'^a_f(\gb) N^{-1/2} \Theta_{a}(x,y) dy $
  \end{enumerate}
  where $C'^a(\gb)$ and $C'^a_f(\gb)$ are positive constants depending on $\gb$. 

  \end{proposition}
  
   \textit{Proof}
    The techniques developped in \cite[Section 7.2]{CD1} can be immediately 
    adapted to our context. Note that a crucial point in this procedure is to use 
    the asymptotics of 
    Lemma \ref{Pr} and the uniform convergence part of it. Hence we just give a very rough sketch of the proof and refer to  \cite[Section 7.2]{CD1}
     for further details. The authors of \cite{CD1} introduced a kernel $A_{x,dy}(n)$ satisfying the following assumptions:
    \begin{enumerate}
  \item the spectral radius of $G_{x,dy} := \sum_{n \in \N} A_{x,dy}(n)$ is strictly smaller than one;
    
  \item  as $n \to \infty$, 
    \begin{equation}\label{kerA}
      A_{x,dy}(n) \sim  L_{x,dy}/n^{2};
    \end{equation}
    furthermore, 
    there exists a positive constant $\cC$ such that for every $x \in [0,a]$ and every  closed set $F \subset [0,a]$, 
    $A_{x,F}(n) \leq  \cC L_{x,F}/n^{2}$;
  \item there exists $ \gga > 1$ such that $((1-\gga G)^{-1} \circ L \circ (1-\gga G)^{-1})_{x,F} < \infty$ (recall that the notation $\circ$ was 
   introduced in Section \ref{ShFT})
  for all $x \in [0,a]$ and for all $F$ Borel subset of $[0,a]$.
  \end{enumerate}
   
  Then, as $n \to \infty$, they proved that the following equivalence holds for $x \in [0,a]$ and $F$ any Borel subset of $[0,a]$:
  \begin{equation}\label{eqk}
    A_{x,F}^{*k}(n) \sim \frac{\left((1-\gga G)^{-1} \circ L \circ (1-\gga G)^{-1}\right)_{x,F}}{n^{2}}.
  \end{equation}
  
   A similar statement can be proved replacing $n^{2}$ by $n^{3/2}$ in both \eqref{eqk} and \eqref{kerA}. The authors of \cite{CD1} apply 
    this statement to the kernel $e^{\gb} F_{x,dy}(n)$ (with the notations of the current paper); all three assumptions above 
     can easily be checked for this kernel also in our context; recall in particular \eqref{defker} and Proposition
      \ref{ExFree} for the first point and Lemma \ref{Pr} for the second point. The verification of the last point relies in their case also 
       on a local limit estimate, which is Lemma \ref{Pr} in our setup,  and we use 
        the fact that the series  $\sum_{n \geq 1} n^{-3/2}$ is convergent
        (instead of the series $\sum_{m \geq 0} \sum_{n \geq m} \frac{1}{mn^{3/2}}$).
     $\qed$

     As it was done in the standard homogeneous case (see \cite[Theorem 2]{CGZ}), making use of Proposition \ref{esti}, we can
     describe the set of contact points in the subcritical regime. Namely, 
  we introduce a probability law $ p_{\gb,N}^f(\cdot, \cdot) $ on  $\{1,\ldots,N\} \times (\R^+)^N$ defined by:
    \begin{equation}
      p_{\gb,N}^f (A,dx) := \frac{1}{Z_{N,a,\gb}^f } 
      e^{ \gb |A|} \prod_{j=1}^{|A|} F_{x_{t_{j-1}},dx_{t_j}}(t_j - t_{j-1}) \ind_{ x_{t_j} \in [0,a], \forall j \in \{ 0, \ldots, |A| \} }
  \end{equation}
    where $t_0 := 0$, $x_0 := 0$ and $A := \{ t_1 < t_2 < \ldots < t_{|A|} \}$.
     
  This law is related to $\gP_{N,a,\gb}^f$ in the following way. We can write
    \begin{equation} \label{nine}
	\gP_{N,a,\gb}^f (dx) = \sum_{ A \subset \{1,\ldots,N\}} \int_{[0,a]^{|A|}} p_{\gb,N}^f (A,dy) \mathbb{P}_{A,y}(dx) 
    \end{equation}
  where $\mathbb{P}^{f}_{A,y}( \cdot )$ is the law of $(S_1,\ldots, S_N)$ conditioned on the event $\cE^{f}_{N,A,y}$ which is defined by:
  \begin{equation}\label{epscond}
    \cE^{f}_{N,A,y} :=  \Big\{ S_i = y_i;i \in A  \cup \{ 0 \} \Big\} \cap \Big\{  S_i > a, i \notin A \Big\}.
  \end{equation}
  
  We define the analogous quantities in the constrained case. Namely, we consider the 
  probability law $ p_{\gb,N}^c(\cdot, \cdot) $
   on $\{1,\ldots,N-1\} \times (\R^+)^N$, which is defined by:
    \begin{equation}\label{ten}
      p_{\gb,N}^c (A,dx) := \frac{1}{Z_{N,a,\gb}^c} e^{ \gb (|A|+1)}
  \prod_{j=1}^{|A|+1} F_{x_{t_{j-1}},dx_{t_j}}(t_j - t_{j-1}) \ind_{ x_{t_j} \in [0,a], \forall j \in \{ 0, \ldots, |A| \} }
  \end{equation}
  where $t_0 := 0$, $t_{|A|+1} := N$, $x_0 := 0$ and $A := \{ t_1 < t_2 < \ldots < t_{|A|} \}$.
   One readily realizes that 
   \begin{equation}\label{elven}
	\gP_{N,a,\gb}^c (dx) = \sum_{ A \subset \{1,\ldots,N\}} \int_{[0,a]^{|A| +1 }} p_{\gb,N}^c (A,dy) \mathbb{P}^c_{A,y}(dx) 
    \end{equation}
  where for $y \in (\R^+)^N$, $\mathbb{P}^{c}_{A,y}( \cdot )$ is the law of $(S_1,\ldots, S_N)$ conditioned
  on the event $\cE^{c}_{N,A,y}$ which is defined by:
  \begin{equation}
    \cE^{c}_{N,A,y} :=  \Big\{ S_i = y_i;i \in A  \cup \{ 0 \} \cup \{ N \} \Big\} \cap \Big\{  S_i > a, i \notin A \cup \{ N \} \Big\}.
  \end{equation}

    For $A \subset \{1,\ldots,N\}$, we define
    \begin{equation}\label{lftarig}
     L(A) = \max(A \cap [0,N/2]) \hspace{6 pt } \text{and} \hspace{6 pt }  R(A) := \min((A \cap [N/2, N]) \cup \{N\}).
    \end{equation} 
    
    Combining the following result with \eqref{nine} and \eqref{elven} implies Theorem \ref{TRLO}: 

    \begin{lemma} \label{MLL}
      For $\gb < \gb_c^a$, the following estimate holds:
      \begin{equation}\label{end}
      \lim_{ L \to \infty } \limsup_{ N \to \infty }  \sup_{ x \in \R^N} p_{\gb,N}^f ( \max A \geq L,dx) = 0.
  \end{equation}
   The corresponding estimates in the constrained case read:
     \begin{equation}
      \lim_{ L \to \infty } \limsup_{ N \to \infty }  \sup_{ x \in \R^N} p_{\gb,N}^c ( L(A) \geq L,dx) = 0
  \end{equation}
  and 
  \begin{equation}\label{env}
    \lim_{ L \to \infty } \limsup_{ N \to \infty }  \sup_{ x \in \R^N} p_{\gb,N}^c ( R(A) \leq N - L,dx) = 0.
  \end{equation}
    \end{lemma}
    
   The proof of these convergences follows by making use of the equivalences from Proposition \ref{esti} and goes along the same lines 
    as the proof of \cite[Propositions 5 and 6]{DGZ}.

    \subsection{Scaling limits in the subcritical regime. The free case.} \label{IVLSR}

  The goal of this Section is to prove Theorem \ref{MAIN} in the free case for $\gb < \gb_c^a$. We stress that similar
  ideas to the ones developed in this Section 
   will be used in the constrained case also (see Section \ref{casconst}), but that this later case is technically more involved.
  
    In what follows, we define $\tau_{(- \infty,0)} := \inf \{ j \geq 0, S_j < 0 \}$. 
     
     Combining the estimates on the contact set of Theorem \ref{TRLO} and 
  the representations of \eqref{nine} and \eqref{ten}, we can restrict the analysis to the trajectories whose 
  contacts with the strip are close to $\{ 0 \}$. After integrating
  over the first step after the last contact with the strip and making use of Markov's property, the remaining process
  is simply the random walk conditioned to stay above the strip. Finally, the convergence towards the brownian meander of Theorem \ref{MAIN} is 
   a consequence of the following result which is due to Shimura:  
  
  \begin{theorem} [Example 4.1 in \cite{Shi}] \label{ShiTT} 
  Let $x_N$ be a positive sequence such that $x_N N^{-1/2}\to0$ as $N \to \infty$. One has the following functional limit convergence:
    \begin{equation}
    \gP_{x_N} \Big[ \hspace{.1cm} \cdot \hspace{.1cm}  \Big| \tau_{(- \infty,0)} > N \Big] \circ \Big(X^N\Big)^{-1} \Longrightarrow m(\cdot).
    \end{equation} 
  \end{theorem}
  For clarity, we summarize the steps  of the proof of Theorem \ref{MAIN} in the next key lemma; then we show 
  that we may apply Lemma \ref{KEYLEm} to our setup, and finally we go to its proof in Section \ref{preuv}.  
  
  \begin{lemma} \label{KEYLEm}
    Let $L$ be a positive integer. Recall that $S^{N}$, the image of $S$ through the application $X^{N}$, was defined in \eqref{def3S}.
  Assume the following assumptions hold:
    \begin{enumerate}
    \item  for any $\gep > 0$, one has
    \begin{equation}
      \lim_{N \to \infty} \mathbf{P}_{N,a,\gb}^f \left[ \sup_{ t \in \left[0,\frac{\max \cA}{N}\right]}S^{N}_t \geq \gep \right] = 0.
    \end{equation}
    \item   for every sequence of subsets $A_{k_{N}} \subset \{1,\ldots,N\}$ such that $\max A_{k_{N}}$ is fixed and satisfies $\max A_{k_{N}} = L$ and 
    for every sequence of vectors $x_{k_{N}} = (x_{k_{N}}^{1},\ldots,x_{k_{N}}^{N}) \in \R^N$, if $S$ follows the law $\mathbb{P}^{f}_{A_{k_{N}},x_{k_{N}}}$,
  one has the convergence in law: 
    \begin{equation}
      \left(\sqrt{1 -\frac{L}{N}} S^N_{\frac{L}{N} + t(1 - \frac{L}{N})} \right)_{ t \in [0,1]} \Rightarrow m
    \end{equation}
      where $m$ denotes the law of the brownian meander. 
    \end{enumerate}
    
      Then one has the weak convergence 
    \begin{equation}
      Q_{N,a,\gb}^f \Rightarrow m.
    \end{equation}
  
  \end{lemma}

  \textbf{ Proof of Theorem \ref{MAIN}.} 
   
    \textit{Verification of the first point of Lemma \ref{KEYLEm}.}
    We write :
    \begin{equation}
  \mathbf{P}_{N,a,\gb}^f \left[ \sup_{ t \in \left[0,\frac{\max \cA}{N}\right]}S^{N}_t \geq \gep \right] 
  = \mathbf{P}_{N,a,\gb}^f \left[ \max_{j=1, \ldots, \max \cA} S_j \geq \gep \gs \sqrt{N}; \max \cA > L \right]  
    \end{equation}
    \begin{equation*}
      + \mathbf{P}_{N,a,\gb}^f \left[ \max_{j=1, \ldots, \max \cA} S_j \geq \gep \gs \sqrt{N}; \max \cA \leq L \right].
    \end{equation*}
     
     Thanks to Theorem \ref{TRLO}, for any fixed $\eta > 0$, one can choose $L_{0} > 0$ such that  
   for every $L \ge L_0$: 
    \begin{equation}\label{primo}
      \limsup_{N \to \infty} \mathbf{P}_{N,a,\gb}^f \left[ \max_{j=1, \ldots, \max \cA} S_j \geq \gep \gs \sqrt{N}; \max \cA > L \right]  \leq \eta/2.
    \end{equation}
    
      Then we note that: 
    \begin{equation}
    \mathbf{P}_{N,a,\gb}^f \left[ \max_{j=1, \ldots, \max \cA} S_j \geq \gep \gs \sqrt{N}; \max \cA \leq L_{0} \right]
  = \frac{\mathbf{E} \left[ \ind_{ \left\{ \max_{j=1, \ldots, \max \cA} S_j \geq \gep \gs \sqrt{N} \right\}} e^{\gb \sum_{i=1}^N \ind_{ S_i \in [0,a]}}
  \ind_{ \max \cA \leq L_{0}} \ind_{ T_1^- > N} \right]}{{Z_{N,a,\gb}^f}}
  \end{equation}
  so that using the estimates on $Z_{N,a,\gb}^f$ from Proposition \ref{esti} and Lemma \ref{fluctu}, we get that there exists a constant $\cC > 0$
 such that : 
    \begin{equation}
    \begin{aligned}
     \mathbf{P}_{N,a,\gb}^f \left[ \max_{j=1, \ldots, \max \cA} S_j \geq \gep \gs \sqrt{N}; \max \cA \leq L_{0} \right] & \leq \cC N^{1/2}e^{\gb L_{0}}
     \mathbf{E} \left[ \ind_{ \max_{j=1, \ldots, L} S_j \geq \gep \gs \sqrt{N}}  \ind_{l_N \leq L_{0}} \ind_{ T_1^- > N} \right] \\
	&  \leq \cC e^{\gb L_{0}} \sum_{j=1}^{L_{0}} N^{1/2} \mathbf{P} \left[ S_j \geq \gep \gs \sqrt{N} \right].
    \end{aligned}
  \end{equation}
  
   Since each term in the sum of the right hand side of the above inequality vanishes in the asymptotic $N \to \infty$, we deduce that 
    \begin{equation}\label{secondo}
      \mathbf{P}_{N,a,\gb}^f \left[ \max_{j=1, \ldots, \max \cA} S_j \geq \gep \gs \sqrt{N}; \max \cA \leq L_{0} \right] \leq \eta/2
    \end{equation} 
  as soon as $N$ is large enough. 
  
   Combining this last inequality with \eqref{primo}, we deduce the first point of Lemma \ref{KEYLEm}. 
    
    $ \qed $
  
  \textit{Verification of the second point of Lemma \ref{KEYLEm}.} 
  
   We first verify it in the case where the sequence $(A_{k_{N}})_{N}$ is constant and satisfies 
   $A_{k_{N}} = \emptyset$ (so that we necessarily have $n_{0} = 1$, and it is easily checked that 
    we do not use sets of trajectories with probability $0$ in the rest of the proof). Let $\gep > 0$.
  We consider a Lipschitz bounded functional $\Phi$ on $\cC([0,1],\R)$, that is such that there
  exist two positive constants $c_1$ and $c_2$ verifying that for every $f,g \in \cC([0,1],\R)$, one has:
  \begin{equation} \label{Ph}
      |\Phi(f)| < c_1 \hspace{.4 cm} \text{and} \hspace{.4 cm} |\Phi(f) - \Phi(g) | \leq c_2 ||f-g||_{\infty}.
  \end{equation}
    Here, the event $\cE^{f}_{N,A_{k_{N}},x_{k_{N}}}$ appearing in \eqref{epscond} is the event $\{ S_1 > a, \ldots, S_N > a \}$;
    conditioning on $S_1$ and using Markov's property, one gets:
  \begin{equation} \label{I} \begin{split}
    &  \mathbb{E}_{A}^f \Big[ \Phi \left( (X^N_{t})_{t \in [0,1]} \right) \Big] \\
  & \phantom{iiii}  = \int_a^{\infty} \frac{\mathbf{E}\Big[ \Phi \left( X^N(t,S_2,\ldots,S_N)
  \right ), S_1 = t, S_2 > a, \ldots, S_N > a  \Big]}{ \gP[S_{1} > a, \ldots, S_{N} > a]} dt.
   \end{split}
  \end{equation}

  Then we use the Markov property and the invariance by translation of $S$ to get that for any $t \geq a$: 
  \begin{equation} 
    \begin{split}
      &   \mathbf{E}\Big[ \Phi \left( X^N(t,S_2,\ldots,S_N)\right ), S_1 = t, S_2 > a, \ldots, S_N > a  \Big] \\
    & \phantom{iiii} = h(t) \mathbf{E}_{t-a} \Big[ \Phi \left( X^N(t,S_1+a, \ldots, S_{N-1}+a) \right), \tau_{(- \infty,0)} > N-1 \Big].
      \end{split}
  \end{equation}

    For any $(x_1, \ldots, x_{N-1}) \in (\R^+)^{N-1}$ and $t \in  [a,N^{1/4}]$, one has 
  \begin{equation}\label{Ko}
    \begin{split}
      & \Big| \Phi \left( X^N(t,x_1+a, \ldots, x_{N-1}+a)\right) - \Phi \left( X^{N-1}(x_1, \ldots, x_{N-1})\right) \Big| \\
  & \phantom{iiiiii}  \leq \frac{c_{2}}{\sqrt{N}} \left( \sup_{j=1, \ldots, N-1} |x_j - x_{j-1}| + a+N^{1/4}\right).
    \end{split}
  \end{equation}
  
    Theorem \ref{ShiTT} implies that for $t \in (a,N^{1/4}), \mathbf{E}_{t-a} \Big[ \hspace{.1 cm}  \cdot \hspace{.1 cm} \Big|
  \tau_{(- \infty,0)} > N-1 \Big] \circ \Big(X^{N-1} \Big)^{-1}$ 
  converges towards $m(\cdot)$.
  In particular, using the tightness criterion of Kolmogorov, this implies  the fact 
  that $  \frac{1}{\sqrt{N}} \left(  \sup_{j=1, \ldots, N-1} |S_j - S_{j-1}| +  a+N^{1/4}\right) =: \mathcal{Y}^a_N$
  converges towards zero in probability when $(S_j)_{j \leq N}$ is distributed according 
  to $ \mathbf{E}_{t-a} \Big[ \hspace{.1 cm}  \cdot \hspace{.1 cm} \Big|
  \tau_{(- \infty,0)} > N-1 \Big]$. 

  Thus, one has:
  \begin{equation}\label{Pi}
  \begin{aligned}
   & \Big|\mathbf{E}_{t-a}\Big[ \Phi \left( X^N(t,S_1+a, \ldots, S_{N-1}+a)\right) \Big| \tau_{(- \infty,0)} > N-1\Big] \\
      & \phantom{iiiiiiiiiiiiiiiiiiiiiii} - \mathbf{E}_{t-a}\Big[ \Phi \left(X^{N-1}(S_1, \ldots, S_{N-1})\right) \Big|
  \tau_{(- \infty,0)} > N-1\Big] \Big| \\
    & \leq \mathbf{E}_{t-a} \Big[ \ind_{ \mathcal{Y}^a_N > \gep} |\Phi \left( X^N(t,S_1+a, \ldots, S_{N-1}+a)\right) \\
  &  \phantom{iiiiiiiiiiiiiiiiiiiiiii} - \Phi \left(X^{N-1}(S_1, \ldots, S_{N-1})\right) | 
   \Big| \tau_{(- \infty,0)} > N-1 \Big] \\
    & + \mathbf{E}_{t-a} \Big[ \ind_{ \mathcal{Y}^a_N \leq \gep} |\Phi \left( X^N(t,S_1+a, \ldots, S_{N-1}+a)\right) \\
  & \phantom{iiiiiiiiiiiiiiiiiiiiiii} - \Phi \left(X^{N-1}(S_1, \ldots, S_{N-1})\right) |
     \Big| \tau_{(- \infty,0)} > N-1 \Big] \\
  & \phantom{iiiiiiiiiiiiiiiiiiiiiii} \leq 2c_1 \mathbf{P}_{t-a} \Big[ \mathcal{Y}^a_N > \gep \Big| \tau_{(- \infty,0)} > N-1 \Big] + c_2 \gep
\end{aligned}
\end{equation}
  where in the last inequality we made use of \eqref{Ko}.
  We finally choose $N$ large enough such that the last term above is smaller than
  say $2 c_{2} \gep$. Informally stated, \eqref{Pi} implies that in the following,
  one can approximate $\mathbf{E}_{t-a}\Big[ \Phi \left( X^N(t,S_1+a, \ldots,S_{N-1}+a)\right) \Big]$
  by $\mathbf{E}_{t-a}\Big[ \Phi \left( X^{N-1}(S_1, \ldots,S_{N-1}) \right) \Big]$.  
  
   We then rewrite \eqref{I} as 
    \begin{equation}\label{abat}
    \begin{aligned}
        &  \mathbb{E}_{A}^f \Big[ \Phi \left( (X^N_{t})_{t \in [0,1]} \right) \Big] \\
  & \phantom{iiii}  = \int_a^{N^{1/4}} h(t)\frac{\gP_{t-a}[\tau_{(- \infty,0)} > N-1]}{\gP[S_{1} > a, \ldots, S_{N} > a]}
  \mathbf{E}_{t-a} \Big[ \Phi \left( X^N(t,S_1+a, \ldots, S_{N-1}+a) \right)\Big| \tau_{(- \infty,0)} > N-1 \Big] dt   \\ 
   & \phantom{iiii} + \int_{N^{1/4}}^{\infty}  h(t)\frac{\gP_{t-a}[\tau_{(- \infty,0)} > N-1]}{\gP[S_{1} > a, \ldots, S_{N} > a]}
   \mathbf{E}_{t-a} \Big[ \Phi \left( X^N(t,S_1+a, \ldots, S_{N-1}+a) \right)\Big| \tau_{(- \infty,0)} > N-1 \Big] dt.
    \end{aligned}
    \end{equation}
    
    For the first term in the right hand side of the above equality, we first 
    replace the expectation term in the integral 
    by $\mathbf{E}_{t-a}\Big[ \Phi \left( X^{N-1}(S_1, \ldots,S_{N-1}) \right) \Big]$, losing a constant $2c_{2}\gep$ by 
    doing so; given the range of integration, 
     by  Theorem \ref{ShiTT}, this last term converges towards $m(\Phi)$. Then we combine the 
  dominated convergence theorem and the fact that
  $\int_a^{\infty} h(t)\frac{\gP_{t-a}[\tau_{(- \infty,0)} > N-1]}{\gP[S_{1} > a, \ldots, S_{N} > a]} dt = 1$ for every $N$ to get 
  that the first term in the right hand side of the above equality converges  
 as $N \to \infty$ towards $m(\Phi)$.
 
  On the other hand, since for $ t \in [N^{1/4},\infty)$, we have $t^{2}/\sqrt{N} \geq 1$, we get that 
 the second term in the right hand side of \eqref{abat} is smaller than 
  \begin{equation}\label{sma}
  c_{2} \int_{N^{1/4}}^{\infty} \frac{ t^{2} h(t)}{N^{1/2} \gP[S_{1} > a, \ldots, S_{N} > a] }dt.
  \end{equation}
  
   By Lemma \ref{stayabove}, the sequence $(N^{1/2} \gP[S_{1} > a, \ldots, S_{N} > a])_{N} $ converges towards a positive limit as $N \to \infty$. 
        
  We make use of this convergence and of the fact that  $\E[X^{2}] < \infty$ to get that the term 
  in \eqref{sma} vanishes as $N \to \infty$; hence the case where the sequence $A_{k_{N}}$ is identically equal to $\emptyset$ is resolved.

  For a generic sequence $(A_{k_{N}},x_{k_{N}})_{N}$, we make use of the Markov property 
 and of what we just proved. More precisely, since for any $N$, we have $A_{k_{N}} \cap [0,L] = A_{k_{N}}$, using Markov's property, for every measurable 
 function $H: \R^{N-L} \to \R$, we get the equality
 \begin{equation} \begin{split}
    \gE\left [ H\left( (S_{L},\ldots, S_{N}) \right ),\cE^{f}_{N,A_{k_{N}},x_{k_{N}}}\right ] & 
    = \gE\left [  H\left( (S_{L},\ldots, S_{N}) \right ), \cE^{f}_{L,A_{k_{N}}\cap[0,L],x_{k_{N}}} \cap \{S_{L+1} > a, \ldots, S_{N} > a \}\right ] \\
    & = \gP\left [ \cE^{f}_{L,A_{k_{N}}\cap[0,L],x_{k_{N}}}\right ] \gE_{x^{|A|}_{k_{N}}}\left [H(S_{0}, \ldots, S_{N-L}) \ind_{S_{1} > a, \ldots, S_{N-L} > a} \right ]. 
     \end{split}
    \end{equation}   
    
  Note that both sides of this equality might be zero in the case where $\gP_{x^{|A|}_{k_{N}}}[\{S_{1} > a, \ldots, S_{N-L} > a\}]=0$.

   From this we deduce:
  \begin{multline}
    \mathbb{E}^{f}_{A_{k_{N}},x_{k_{N}}} \left[ \Phi \left( \left(\sqrt{1 -\frac{L}{N}} S^N_{\frac{L}{N} + t(1 - \frac{L}{N})} \right)_{ t \in [0,1]} \right) \right] \\
    =  \mathbb{E}^{f}_{x^{|A|}_{k_{N}}} \Big[ \Phi \left((S_t^{N-L})_{t \in [0,1]}\right) \Big| S_1 > a, \ldots, S_{N-L} > a \Big]. 
  \end{multline}
 Finally, we note that for all $x \in [0,a]$, one has the equality:
    \begin{equation}
      \mathbf{E}_{x} \Big[ \Phi \left( S^{N-L}_{t})_{ t \in [0,1]} \right) \Big| S_1 > a, \ldots, S_{N-L} > a \Big] = \mathbb{E}^{f}_{ \emptyset, x} \Big[ \Phi \left( S^{N-L}_{t})_{ t \in [0,1]} \right) \Big].  
    \end{equation}
    
     We already proved that the right hand side in the above equality converges towards $m(\Phi)$ in the particular case $x = 0$. Getting the 
     same convergence for any $x \in [0,a]$ works in the same way, and hence we get the second point of Lemma \ref{KEYLEm}. $ \qed $

  \subsubsection{Proof of Lemma \ref{KEYLEm}} \label{preuv}
    We consider $\gep, \eta >0$, $L_0$ a positive integer and $\Phi$ a continuous function on $C([0,1],\R)$. We write:
    \begin{equation}\label{Fin}
    \begin{aligned}
       Q_{N,a,\gb}^f \Big[ \Phi \Big] 
     = \sum_{l=0}^{L_0} \sum_{ A \subset \{1,\ldots,N\}; \max A = l} \int_{[0,a]^{|A|}}
  p_{\gb,N}^f \left(A,dx\right) \mathbb{P}^{f}_{A,x} \Big[ \Phi(S^N) \Big] + Q_{N,a,\gb}^f \Big[ \Phi \ind_{ \max \cA > L_0} \Big].
    \end{aligned}
  \end{equation}
  
   We first prove that, for all $ A \subset \{1,\ldots,N\}$ such that $\max A \leq L_{0}$, we have the convergence
   \begin{equation}
    \mathbb{P}^{f}_{A,x} \Big[ \Phi(S^N) \Big] \to m(\Phi).
   \end{equation} 
    
  We note $L$ for the quantity $\max A$ and for notational convenience we
  write $f_N(t) := L/N + t(1-L/N)$ and $g_N(t) := \frac{(t-L/N)}{1-L/N}$ its
  inverse (and we set $f_0(t) = g_0(t) = t$). 

    We first note that for every $ n > 0$, for every $ t_1 < t_2 < \ldots < t_n \in [0,1]^n$ and for every continuous bounded function
  $F : [0,1]^n \to \R$, one has the convergence 
    \begin{equation}
  \mathbb{P}^{f}_{A,x}\left[F(S^{N}_{t_{1}},S^{N}_{t_2}, \ldots, S^{N}_{t_n}) \right] \to m\left[F( \go_{t_1},\ldots, \go_{t_n})\right],
    \end{equation}
     where $(\go_{t})_{t \in [0,1]}$ denotes the canonical process under the law $m$.
      
  Indeed, since $F$ is continuous and bounded, by dominated convergence, as $N \to \infty$, we get:
    \begin{multline}
      \Bigg|\mathbb{P}^{f}_{A,x} \Big[F\left(S^{N}_{t_{1}},S^{N}_{t_2}, \ldots, S^{N}_{t_n}\right) \Big] \\
  - \mathbb{P}^{f}_{A,x}\left[ F\left( \sqrt{1 -\frac{L}{N}} S^N_{f_N(t_1)},
  \ldots, \sqrt{1 -\frac{L}{N}} S^N_{f_N(t_n)}\right)  \right] \Bigg| \to 0.
    \end{multline}

     Since the convergence of
  the second term above towards $m(F(\go_{t_1}, \ldots, \go_{t_n}))$ is the hypothesis 2 of Lemma \ref{KEYLEm}, the
  finite dimensional convergence is proven.
  
    We are left with proving the tightness of the sequence $S^N$ under the law $\mathbb{P}^{f}_{A,x}$, for $A \subset \{1,\ldots,N\}$ such that
  $\max A \leq L_0$. For this, for $\delta > 0$ and for a continuous function $f$ on $[0,1] \to \R^+$ verifying $\sup_{t \in [0,\delta]} f(t) \leq \gep$,
  we introduce its $\delta$-cut counterpart $f^{(\delta)}$; namely, $f^{(\delta)} (x) = \frac{x f(\delta)}{\delta} \ind_{x \in [0,\delta]} + f(x)\ind_{x \geq \delta}$.
  Clearly, we have $ ||f^{(\delta)} - f ||_{\infty} \leq \gep$.

    We combine the $\ga$ H\"older regularity  of 
  the brownian motion for any $\ga \in (0,1/2)$ (see for example \cite[Corollary 1.20]{MP} for a proof of this classical result) and the representation (see 
  for example \cite{DIM})
   \begin{equation}
    (m(t))_{t \in [0,1]} \stackrel{\cL}{=} \left(\frac{1}{\sqrt{1-\kappa_{1}}} |B_{\kappa_{1} + t(1-\kappa_{1})}|\right)_{t \in [0,1]}
   \end{equation} 
   where $\kappa_{1} = \sup_{s \leq 1} \{B_{s} =0 \}$  to get that, for $C$ large enough (recall that $\kappa_{1}$ follows the arcsine law, and in particular
    $\gP[ \kappa_{1} > 1-\eta ]$ can be made arbitrarily small by choosing $\eta$ small enough),
  one has $m( \mathcal{B}_{C}) \geq 1-\gep $ where 
  \begin{equation}
    \mathcal{B}_C := \Big\{ f \in C([0,1],\R), \sup_{x,y \in [0,1]} \frac{|f(x) - f(y)|}{|x-y|^{1/3}}  \leq C  \Big\}. 
  \end{equation}
   
    Therefore for such a $C$ and for $N$ large enough, applying hypothesis 2 of Lemma \ref{KEYLEm}, we get:
  \begin{equation}
    \mathbb{P}^{f}_{A,x} \left[ \left(\sqrt{1 -\frac{L}{N}} S^N_{f_N(t)} \right)_{t \in [0,1]} \in \mathcal{B}_C \right] \geq 1 - 2 \gep.
  \end{equation}
   
    Now we are ready to prove the Kolmogorov criterion for $S^N$ under the law $\mathbb{P}^{f}_{A,x}$. 
    We have to show that for any  given $\gd > 0$, there exists $N_0$ such that:
  \begin{equation} \label{AMPBB}
    \mathbb{P}^{f}_{A,x} \Big[ \sup_{ s,t, |s-t| \leq \gd} |S_s^N - S_t^N| \geq \gep \Big] \leq \eta, \hspace{.1 cm}\text{for all }  N \geq N_0.
  \end{equation}
    
    We claim that we can restrict ourselves to show \eqref{AMPBB} by
    replacing $S^N$ by its $L/N$-cut counterpart, which we denote by  $\tilde{S}^N$. Let us prove this claim. 
    
     Indeed, since for $\kappa \in (0,1)$, $f$ and $f^{(\kappa)}$ coincide on $[\kappa,1]$, we have:
     \begin{equation}
     \begin{aligned}
   \mathbb{P}^{f}_{A,x} \Big[ \sup_{ s,t, |s-t| \leq \gd} |S_s^N - S_t^N| \geq \gep \Big] 
  & = \mathbb{P}^{f}_{A,x} \Big[ \sup_{ s,t \geq L/N, |s-t| \leq \gd} |\tilde{S}_s^N - \tilde{S}_t^N| \geq \gep \Big] \\
   & \phantom{iiiiiii} + 2 \mathbb{P}^{f}_{A,x} \Big[ \sup_{ t \leq L/N \leq s, |s-t| \leq \gd} |S_s^N - S_t^N| \geq \gep \Big] \\ 
   & \phantom{iiiiiii} + \mathbb{P}^{f}_{A,x} \Big[ \sup_{ t \vee s \leq L/N, |s-t| \leq \gd} |S_s^N - S_t^N| \geq \gep \Big]. 
     \end{aligned}
     \end{equation}

    Since
     \begin{equation}
     \begin{aligned}
          \mathbb{P}^{f}_{A,x} \Big[ \sup_{ t \vee s \leq L/N, |s-t| \leq \gd} |S_s^N - S_t^N| \geq \gep \Big]
      \leq 2 \mathbb{P}^{f}_{A,x} \Big[ \sup_{t \leq L/N} S_t^N > \gep/2  \Big],
     \end{aligned}
     \end{equation} 
     we can make use of the first item of Lemma \ref{KEYLEm} to deduce that this term vanishes for $N \to \infty$. By triangular inequality, we also have
      \begin{equation}
      \begin{aligned}
   & \mathbb{P}^{f}_{A,x} \Big[ \sup_{ t \leq L/N \leq s, |s-t| \leq \gd} |S_s^N - S_t^N| \geq \gep \Big] \\
   & \leq \mathbb{P}^{f}_{A,x} \Big[ \sup_{ L/N \leq s \leq L/N + \delta} |S_s^N - S_{L/N}^N| \geq \gep/2 \Big] 
   + 2 \mathbb{P}^{f}_{A,x} \Big[ \sup_{t \leq L/N} S_t^N > \gep/4  \Big],
      \end{aligned}
      \end{equation} 
      which  entails the claim, by using once again the first item of Lemma \ref{KEYLEm} and the obvious inequality 
    \begin{equation}
      \mathbb{P}^{f}_{A,x} \Big[ \sup_{ L/N \leq s \leq L/N + \delta} |S_s^N - S_{L/N}^N| \geq \gep/2 \Big] \leq \mathbb{P}^{f}_{A,x}
      \Big[ \sup_{ s,t \geq L/N, |s-t| \leq \gd} |\tilde{S}_s^N - \tilde{S}_t^N| \geq \gep/2 \Big].
    \end{equation}

    Thus we are left with showing that there exists $\gd > 0$ such that for $N$ large enough, one has:
    \begin{equation} \label{A}
      \mathbb{P}^{f}_{A,x} \Big[ \sup_{ s,t > L/N, |s-t| \leq \gd} |\tilde{S}_s^N - \tilde{S}_t^N| \geq \gep \Big] \leq \eta.
    \end{equation}
     
  Now we write:
  \begin{equation}
      \Big|\tilde{S}_s^N - \tilde{S}_t^N \Big| \leq \Big|\tilde{S}_{f_N(s)}  ^N - \tilde{S}_{f_N(t)}^N\Big| +
  \Big|\tilde{S}_{f_N(s)}  ^N - \tilde{S}_s^N\Big| + \Big|\tilde{S}_{f_N(t)}  ^N - \tilde{S}_t^N\Big| 
  \end{equation} 
    so that, for every $\gd > 0$,
  \begin{equation}
  \begin{split}
    & \mathbb{P}^{f}_{A,x} \Big[ \sup_{ s,t > L/N, |s-t| \leq \gd} \Big| \tilde{S}_s^N - \tilde{S}_t^N| \geq \gep \Big] \\
  & \phantom{iiiii} \leq \mathbb{P}_{A,x} \Big[
  \sup_{  |s-t| \leq \gd} \Big|\tilde{S}_{f_N(s)}  ^N - \tilde{S}_{f_N(t)}^N \Big| \geq \gep/3 \Big] 
      + 2 \mathbb{P}^{f}_{A,x} \Big[ \sup_{  s \in [0,1] } \Big|\tilde{S}_{f_N(s)}  ^N - \tilde{S}_s^N \Big| \geq \gep/3 \Big]. 
    \end{split}
    \end{equation}
     
  The first term in the right hand side of the above inequality can be made smaller than $\eta/2$ for $\gd$ small
  enough as soon as $N$ is large enough using the second hypothesis of Lemma \ref{KEYLEm}. 
  For the second term, we get
  \begin{equation} \label{UNDERNIER}
  \begin{split}
  & \mathbb{P}^{f}_{A,x} \Big[ \sup_{  s \in [0,1] } \Big|\tilde{S}_{f_N(s)}^N - \tilde{S}_s^N\Big| \geq \gep/3 \Big] \\
    & \phantom{iiiiiii} = \mathbb{P}^{f}_{A,x} \Big[ \sup_{  s \in [0,1] }  \Big|\tilde{S}_{f_N(s)}  ^N - \tilde{S}_{s}^N  \Big| \geq \gep/3;
  \left(\tilde{S}_{f_N(t)}^N\right)_{t \in [0,1]} \in \mathcal{B}_C  \Big] \\
    & \phantom{iiiiiiiiiiiiiiiiiii} + \mathbb{P}^{f}_{A,x} \Big[ \sup_{  s \in [0,1] } \Big|\tilde{S}_{f_N(s)}  ^N - \tilde{S}_{s}^N\Big| \geq \gep/3;
  \Big(\tilde{S}_{f_N(t)}^N\Big)_{t \in [0,1]} \notin \mathcal{B}_C \Big].
  \end{split}
  \end{equation}
    
  The last term of equation \eqref{UNDERNIER} above can be made smaller than $\eta/3$ for $N$ large enough since $\mathcal{B}_C$ is an $m$ continuity set
  (that is a set whose boundary is of null $m$ measure) and by using the portemanteau theorem, which states that in this case 
    \begin{equation}
      \mathbb{P}^{f}_{A,x} \Big[ (\tilde{S}_{f_N(t)}^N)_{t \in [0,1]} \in \mathcal{B}_C \Big] \to m( \mathcal{B}_C)
    \end{equation}
  as $N \to \infty$.  

    Finally, for $(\tilde{S}_{f_N(t)}^N)_{t \in [0,1]} \in \mathcal{B}_C$, for any $s \in [0,1]$, we have 
  \begin{equation}
     \Big|\tilde{S}_{f_N(s)}  ^N - \tilde{S}_{s}^N \Big| \leq C  \Big| f_{N}(s) - s\Big|^{1/3}
  \end{equation}
    and $\sup_{s \in [0,1]} \Big| f_{N}(s) - s\Big|^{1/3}\leq (L_{0}/N)^{1/3}$. Thus, as soon as $N$ is large enough, we have:
  
    \begin{equation}
      \mathbb{P}^{f}_{A,x}\Big[ \sup_{  s \in [0,1] } \Big|\tilde{S}_{f_N(s)}  ^N - \tilde{S}_{s}^N\Big| \geq \gep/3;
  \left(\tilde{S}_{f_N(t)}^N\right)_{t \in [0,1]} \in \mathcal{B}_C \Big] =0
    \end{equation}
      which proves \eqref{A}. Thus we have shown that $\mathbb{P}^{f}_{A,x}\Big[ \Phi(S^N) \Big] \to m(\Phi)$.

  Now we make use of equation \eqref{Fin} and the triangle inequality to get that 
  \begin{equation}\label{ddd}
  \begin{split}
    & \Big|Q_{N,a,\gb}^f [ \Phi ] - m(\Phi)\Big| \\
  & \phantom{iii} \leq  \sum_{l=0}^{L_0} \sum_{ A \subset \{1,\ldots,N\}; \max A = l} \int_{[0,a]^{|A|}} p_{\gb,N}^f \left(A,dx\right)
  \Big|\mathbb{P}^{f}_{A,x} \left[ \Phi\left(S^N \right) \right]- m\left(\Phi\right)\Big| \\
    & \phantom{iiiiiiiiii} +  m(|\Phi|)Q_{N,a,\gb}^f \Big[ \ind_{ \max \cA > L_0} \Big]
    + Q_{N,a,\gb}^f \Big[ |\Phi| \ind_{ \max \cA > L_0} \Big], 
  \end{split}
  \end{equation}
    where we also used the equality 
       \begin{equation}
         \sum_{l=0}^{\infty} \sum_{ A \subset \{1,\ldots,N\}; \max A = l} \int_{[0,a]^{|A|}} p_{\gb,N}^f (A,dx) = 1.
         \end{equation} 
          Since $\Phi(\cdot)$ is bounded, combining this last equality with
    the dominated convergence theorem and
  the fact that $\mathbb{P}^{f}_{A,x} \Big[ \Phi(S^N) \Big] \to m(\Phi)$,
  we deduce finally Theorem \ref{MAIN} for the free case by considering $L_0$ large enough and by using Theorem \ref{TRLO}. $ \qed $

  \subsection{Scaling limits in the subcritical regime. The constrained case}\label{casconst}
  
  The strategy in this Section is similar to the one of the preceeding one, and we choose to skip some of the proofs for lightness. We 
   first mention that the analogous of
    Shimura's result has been shown recently for the normalized excursion in \cite[Corollary 2.5]{CaCha}.

  \begin{theorem}  [Corollary 2.5 in \cite{CaCha}] \label{CC}
  Let $x_N$ and $y_N$ two positive sequences such that both $x_N/\sqrt{N}$ and $y_N/\sqrt{N}$ vanish as $N \to \infty$. 
    One has the following weak convergence:
      \begin{equation}
      \gP_{x_N}\Big[ \hspace{.1cm} \cdot \hspace{.1 cm} \Big| S_N = y_N,\tau_{(- \infty,0)} > N \Big]
  \circ \Big(X^N\Big)^{-1}   \Rightarrow e(\cdot).
      \end{equation}
  \end{theorem}

    Like we did in the free case, we first give a technical lemma which implies the convergence in 
  the constrained case of Theorem \ref{MAIN}.
    \begin{lemma} \label{ccas1}
    Let $L < R$ be positive integers and 
      assume that the following hold:
    \begin{enumerate}
    \item  for any $\gep > 0$, one has 
  \begin{equation}
    \lim_{N \to \infty} \gP_{N,a,\gb}^c \left[ \sup_{t \in [0,L(\cA)/N] \cup [R(\cA)/N,1]} S^{N}_t \geq \gep\right] = 0,
  \end{equation}
   where the variables $(L(\cA), R(\cA))$ were defined in
      \eqref{lftarig}.
  \item  for every sequence of subsets $A_{k_{N}} \subset \{1,\ldots,N-1\}$ such that the couple $(L(A_{k_{N}}),R(A_{k_{N}}))$
    is fixed and satisfies $(L(A_{k_{N}}),R(A_{k_{N}})) = (L,R)$ and 
    for every sequence of vectors $x_{k_{N}} = (x_{k_{N}}^{1},\ldots,x_{k_{N}}^{N}) \in \R^N$, if $S$ follows the law $\mathbb{P}^{f}_{A_{k_{N}},x_{k_{N}}}$,
  one has the convergence in law: 
    \begin{equation}
      \left(\sqrt{\frac{R-L}{N}} S^N_{\frac{L}{N} + t\left( \frac{R-L}{N}\right)} \right)_{ t \in [0,1]} \Rightarrow e
    \end{equation}
      where $e$ denotes the law of the normalized brownian excursion. 
    \end{enumerate}
      Then one has the second convergence of Theorem \ref{MAIN}.
    \end{lemma}

     The proof of Lemma \ref{ccas1} closely follows the one of Lemma \ref{KEYLEm}, so that we choose to skip it. 
  
  \subsubsection{ Proof of Theorem \ref{MAIN} in the constrained case.}
  
    We show that the hypothesis of Lemma \ref{ccas1} are fulfilled. 
  
    \textit{Verification of the first point of Lemma \ref{ccas1}.} 
    
    Combining the equivalence: 
    \begin{equation} \label{App}
      \mathbf{P} [S_N \in [0,a]; S_j > 0, j \leq N] \sim \frac{\int_0^aU(u)du}{\sqrt{2 \pi} \gs N^{3/2}}
    \end{equation}
    which follows from \eqref{const}, and the asymptotics on $Z_{N,\gb,a}^c$ in Proposition \ref{esti}, the proof of this point
    goes very much along the same lines 
     as in the constrained case by using standard facts on the normalized excursion instead of the meander, so that once again
     we choose to skip it. Note that, like in the free case, this proof relies heavily on \eqref{CCCCC} in Theorem \ref{TRLO}. 
     
%

    \textit{Verification of the second point of Lemma \ref{ccas1}.} 
    
     Here we make use of Theorem 
  \ref{CC} in a crucial way. We also first treat the case where the sequence $ (A_{k_{N}})_{N}$ is identically equal to $\emptyset$ (hence once again $n_{0} =1$). 
  We consider $\gep > 0$ and $\Phi$ a Lipschitz bounded functional on $\cC([0,1],\R)$ verifying the same properties as
  in \eqref{Ph}. We write:
    \begin{equation}\label{asplit}
    \begin{split} 
    &\mathbb{E}^c_{A_{k_{N}},x_{k_{N}}} \Big[\Phi\left(X_t^N\right)_{t \in [0,1]}\Big] \\
    & \phantom{i} = \int_{  t,t' \in [a,\infty)^2,   u \in [0,a]}
     \frac{\mathbf{E}\Big[ \Phi\left(X^N(t,S_2, \ldots,S_{N-2},t',u )\right),
  S_1 = t, S_2 > a, \ldots, S_{N-2} > a, S_{N-1}=t' ,S_N = u\Big] }{ \gP\Big[ S_1 > a, \ldots, S_{N-1} > a, S_N \in [0,a] \Big] } dt dt' du. 
    \end{split}
  \end{equation}
   
   For $t,t' \geq a, u \in [0,a]$, we use twice Markov's property to get:
  \begin{equation} 
   \begin{aligned}
    & \mathbf{E}\Big[ \Phi\left(X^N(t,S_2, \ldots,S_{N-2},t',u )\right),
  S_1 = t, S_2 > a, \ldots, S_{N-2} > a, S_{N-1}=t' ,S_N = u\Big] \\
    & \phantom{ii}=  h(t) h(u-t')  \mathbf{E}_{t-a}\Big[ \Phi\left(X^N(t,S_1 +a, \ldots,S_{N-2}+a,t',u )\right) | \tau_{(- \infty,0) > N-2}, S_{N-2} = t'\Big] \\
  & \phantom{iiiiii} \times \gP_{t-a}\left [ \tau_{(- \infty,0) > N-2}, S_{N-2} = t' \right ].  
   \end{aligned}
  \end{equation}
  
   As in \eqref{Ko}, for $\gep > 0$ and for any $(x_{1}, \ldots, x_{N-2}) \in (\R^{+})^{N-2}, (t,t') \in (a, \gep \sqrt{N})^{2}, u \in [0,a]$, we have:
     \begin{equation}\label{Ko1}
    \begin{aligned}
      & \Big|  \Phi\left(X^N(t,x_1 +a, \ldots,x_{N-2}+a,t',u)\right) - \Phi \left( X^{N-2}(x_1, \ldots, x_{N-2})\right) \Big| \\
  & \phantom{iiiiii}  \leq \frac{1}{\sqrt{N-2}} \left(  c_2 \sup_{j=1, \ldots, N-2} |x_j - x_{j-1}|  + c_2 \left(a+2 \gep \sqrt{N}\right)\right).
    \end{aligned}
  \end{equation}
  
   Since Theorem \ref{CC} asserts that
   \begin{equation}
  \lim_{\gep \to 0^{+}} \lim_{N \to \infty} \mathbf{E}_{t-a}\Big[ \Phi\left(X^{N-2}(S_1, \ldots,S_{N-2} )\right) | \tau_{(- \infty,0) > N-2}, S_{N-2} = t'\Big] = e( \Phi),
   \end{equation}
    we can deduce from \eqref{Ko1} as in the free case that
   \begin{equation}\label{exc}
   \lim_{\gep \to 0^{+}} \lim_{N \to \infty} \mathbf{E}_{t-a}\Big[ \Phi\left(X^N(t,S_1 +a, \ldots,S_{N-2}+a,t',u )\right) | \tau_{(- \infty,0) > N-2}, S_{N-2} = t'\Big] = e( \Phi).
   \end{equation}
   
    We split the integral over $[a,\infty)^{2}$ appearing in \eqref{asplit}:   
  \begin{equation}
    \mathbb{E}^c_{A,x} \Big[\Phi\left(X_t^N\right)_{t \in [0,1]}\Big] = 
     \int_{u \in [0,a]} \left ( \int_{(t,t') \in \cD_{1}^{N}} \ldots +  \int_{(t,t') \in \cD_{2}^{N}} \ldots  \right ) 
     \end{equation}   
    where, given $\gep > 0$ and $C > 0$, we defined 
     \begin{equation}\begin{split}
    & \cD_{1}^{N} = \left\{ (t,t') \in [a,\gep \sqrt{N}]^{2} \right\}, \\
    & \cD_{2}^{N} =  \left\{ (t,t') \in \R^{+}, t \vee t' \geq \gep \sqrt{N}  \right\}. 
     \end{split} 
      \end{equation} 
     
     As we proceeded in the free case, making use of the equivalence \eqref{const} and of the convergence \eqref{exc}, we deduce that
    \begin{equation} \begin{split}
    & \lim_{\gep \to 0^{+}} \lim_{N \to \infty} \int_{u \in [0,a]} \int_{(t,t') \in \cD_{1}^{N}} h(t) h(u-t') \frac{\gP_{t-a}\left [ \tau_{(- \infty,0) > N-2}, S_{N-2} = t' \right ]}{\gP\Big[ S_1 > a, \ldots, S_{N-1} > a, S_N \in [0,a] \Big]} \\
    & \phantom{iii} \times  \mathbf{E}_{t-a}\Big[ \Phi\left(X^N(t,S_1 +a, \ldots,S_{N-2}+a,t',u )\right) | \tau_{(- \infty,0) > N-2}, S_{N-2} = t'\Big] dt dt' du \\
     & \phantom{iiiiii} = e(\Phi). 
       \end{split}
      \end{equation}  
    
    Since $\Phi$ is bounded, we are left with showing that 
    \begin{equation}\label{lastconv}
      \lim_{\gep \to 0^{+}} \lim_{N \to \infty} \int_{u \in [0,a]} \int_{(t,t') \in \cD_{2}^{N}} h(t) h(u-t') 
      \frac{\gP_{t-a}\left [ \tau_{(- \infty,0) > N-2}, S_{N-2} = t' \right ]}{\gP\Big[ S_1 > a, \ldots, S_{N-1} > a, S_N \in [0,a] \Big]} dt dt' du = 0.
     \end{equation}

     We show the convergence of the term appearing in the integral of 
     \eqref{lastconv} pointwise for $u \in [0,a]$ (and indeed we just show it for $u = 0$), \eqref{lastconv} follows by dominated convergence.
     
     Since $h(\cdot)$ is bounded, by Gnedenko's local limit theorem, we have 
    \begin{equation}
  \sup_{n \in \N}  \sup_{t \in \R} \sqrt{n}\gP[S_{n} = t] < \infty.
    \end{equation} 
    
     We then recall that a consequence of Lemma \ref{Pr} is the fact that
  \begin{equation}\label{cccon}
  N^{3/2} \gP\Big[ S_1 > a, \ldots, S_{N-1} > a, S_N \in [0,a] \Big] \to \int_{[0,a]^{2}} \Theta_{a}(x,y) dx dy  > 0 
  \end{equation}
    as $N \to \infty$.

     We are then left with showing that 
     \begin{equation}\label{d2}
    \lim_{\gep \to 0^{+}} \lim_{N \to \infty} N \int_{\cD_{2}^{N}} h(t) h(t') dt dt' = 0.
       \end{equation} 
       
     Using symmetry, we get:  
      \begin{equation}
      \begin{aligned}
       \int_{\cD_{2}^{N}} h(t) h(t') dt dt' & \leq \frac{1}{\gep^{2}N} \int_{\cD_{2}^{N}} (t \vee t')^{2} h(t) h(t')dt dt'\\
        & \leq \frac{2}{\gep^{2}N} \int_{(u,v) \in (\R^{+})^{2}, u \geq \gep \sqrt{N}, u \geq v} u^{2} h(u) h(v)du dv\\
        & \leq \frac{2}{\gep^{2}N} \int_{\gep \sqrt{N}}^{\infty} u^{2} h(u) du,
      \end{aligned}
      \end{equation}
     and recalling that $\E[X^{2}] < \infty$, we immediately get \eqref{d2}.
     
      To conclude the proof of Theorem \ref{MAIN}, we are left with dealing with the case of a generic sequence  $A_{k_{N}} \subset \{1,\ldots,N\}$, which 
      is done similarly to the free case.

   $ \qed $
   
   \appendix
   
   \section{The Markov renewal theorem on a general state space.}\label{apA}
   
    In this appendix, we prove the key result which has been used in Section \ref{TLP}.
      To get estimates on $Z^c_{N,a,\gb} = \cP^{\gb}[N \in \tau]$ in the localized phase, we need
    to show an analogous to the classical Markov renewal theorem (which can be found for example in \cite{As})
    in the case where 
  the state space of the $J_i$'s is not countable.
    Surprisingly enough, we have not been able to find a proof of 
    such a natural result in the literature, so that we choose to prove it in a more general context for later reference.  
   
    \begin{theorem}\label{MarRen}
      Let $\cK$ be a compact metric space  and 
      consider a  Markov Renewal process $(J,\tau)$ on $\N \times \cK$ with law $\cP$ and invariant measure $\mu$ such that $\mu$ is strictly 
       positive on $\cK$. Denote by  $K_{x,dy}(n)$ its  transition kernel and assume that 
        \begin{equation}
         \int_{\cK} \sum_{n \in \N} K_{x,dy}(n) = 1
        \end{equation} 
         and 
        \begin{equation}
           \varXi := \int_{(x,y) \in \cK^{2}} \mu(dx) \sum_{n \in \N} n K_{x,dy}(n) < \infty.
        \end{equation}
         
         Denoting by $||\cdot||$ the total variation norm on $\cK$, for any initial distribution $\gl$ of $J_{0}$,
         the following convergence holds:
         \begin{equation}
         \lim_{ N \rightarrow \infty} \cP_{\gl} \left[ \exists k \in \N, \tau_k = N, J_k \in dy  \right] = 
  \frac{ \mu(dy)}{  \varXi }.
         \end{equation} 
     \end{theorem}
     
      \textit{Proof}
     
     We prove Theorem \ref{MarRen} by making use of the ergodic properties of the \textit{forward Markov chain} naturally
    linked to the Markov renewal.

    We consider the Markov process $(A_k,J'_k)_{k \geq 0}$ on $\mathbb{N} \times \cK$ such that $J_{0}'$ is distributed 
    according to $\gl$ and whose transitions are given by:
  \begin{equation}
    P_{\gl} \left[ A_j =k ,J'_j \in dy | A_{j-1} = l, J'_{j-1} = x \right] := \gd_{k,l-1} \gd_{x} (dy) 
  \end{equation}
    if $l \geq 2$ (where $\gd_{x}(\cdot)$ is the Dirac measure concentrated on $\{x\}$) and by 
  \begin{equation}
    P_{\gl} \left[ A_j =k ,J'_j \in dy | A_{j-1} = 1, J'_{j-1} = x \right] := K_{x,dy}(k) \hspace{6 pt}\text{otherwise}.
  \end{equation}
  
  Note that this Markov chain is nothing but the well known \textit{forward recurrence chain} associated to the Markov renewal. 
  In words, $A_{i}$ denotes the time one has to wait from time $i$ until the next renewal happens
  (that is $A_{i} = \inf \{ k > i, \exists j, \tau_{j} = k \} - i$), the Markov chain $J'$ keeping track of
  the last location of its modulating chain.  

    We introduce the probability measure on $\mathbb{N} \times \cK$ defined  by : 
  \begin{equation}
  \Pi(i,dy) := \frac{1}{\varXi } \int_{x \in \cK} \mu(dx) \sum_{j \geq i} K_{x,dy}(j) .
  \end{equation}

    One readily realizes that $\Pi(\cdot,\cdot)$ is the invariant probability of the Markov process $(A_k,J'_k)_{k \geq 0}$. 
  Indeed, for all $(i,y) \in \mathbb{N} \times \cK$, we check:
  \begin{equation}\label{OEIY}
  \Pi P(i,dy) = \Pi(i+1,dy) + \int_{x \in \cK} \Pi(1,dx) K_{x,dy}(i) =  \Pi(i,dy),
  \end{equation}
  where in the second equality we used the fact that  $\mu(\cdot)$ is the invariant probability for 
    the Markov process $(J_k)_{k \geq 1}$ (noting that  $\Pi(1,dx)$ is a multiple of $\mu(dx)$, 
  this is exactly saying that $  \int_{x \in \cK} \mu(dx) K_{x,dy}(i) = \mu(dy) $, 
  which is the second part of equation \eqref{OEIY}). 

  Making use of the positivity of $\mu$ on $\cK$ and of the compactness of $\cK$,
   the Markov chain  $(A,J')$ satisfies the hypothesis of the classical ergodic Theorem (see for example \cite[Theorem 13.3.3]{MT}), so that 
    \begin{equation}
    \left| \left| P_{\gl}^n - \Pi \right | \right| \rightarrow 0
  \end{equation}
  as $ n \rightarrow \infty$. 
  This implies that, as $j \to \infty$, the following convergence holds in total variation norm :
  \begin{equation} \label{nu}
    P_{\gl} \left[ A_j =1 ,J'_j \in dy  \right] \rightarrow \frac{ \mu(dy)}{ \varXi } 
  \end{equation}
  and since $ P_{\gl} \left[ A_j =1 ,J'_j \in dx \right] = \cP_{\gl} \left[ \exists k \in \N, \tau_k =j, J_k \in dx \right] $,
  the proof of Theorem \ref{MarRen} is complete.
  $ \qed $
   
  \section{ The particular case of the $(p,q)$ random walk. }\label{valcr}

     In this appendix, we illustrate how our results can be used in some particular cases linked to 
     the discrete $(p,q)$ random walk $S$. We first stress that the techniques developed in this paper can all be applied in the 
       discrete setup (the adaptations being straightforward), and for example the scaling limits of Theorem \ref{MAIN} hold. 
       
        We consider a $(p,q)$ random walk $S$, that 
       is a symmetric random walk with increments in $\{-1,0,1\}$ such that
      $\gP[S_{1} = 1] = p = \frac{1-\gP[S_{1}=0]}{2} = \frac{1-q}{2} $ with $p \in (0,1/2)$, and $a \in \Z^{+}$.
      Considering the kernel $K(j) = \gP[S_{1} > 0, \ldots, S_{j-1} > 0, S_{j} = 0]$, it can be deduced from 
       \cite[Proposition A.10]{GB} that
       \begin{equation}\label{cstKpq}
        K(j) \sim \sqrt{\frac{p}{8 \pi}} j^{-3/2} \sim c_{K} j^{-3/2} , j \to \infty
       \end{equation} 
         and an easy computation yields 
          \begin{equation}
           \sum_{j \geq 1} K(j) = \frac{1+q}{2}. 
          \end{equation} 
           
      In the 
      homogeneous pinning case, that is $a = 0$,
      it is shown in \cite[Chapter 2]{GB}, that $\gb_{c}^{0} = -\log(1-p)$.
      
      We show here the following result:
      
       \begin{theorem}\label{caspq}
      For the $(p,q)$ random walk, the following equalities hold:
      \begin{equation}\label{pq1}
     \gb_{c}^{1} = -\log\left(1 - \frac{3-\sqrt{5}}{2}p\right)
    \end{equation} 
     and 
      \begin{equation}\label{pq2}
     \gb_{c}^{2} = -\log\left(1-\frac{5-r}{3}p\right)
    \end{equation} 
     where $r = 2\sqrt{7} \cos\left(\frac{\pi}{3} - \frac{\arctan(3\sqrt{3})}{3}\right) \approx 4.405812$. 
     
     Moreover, in the case $a=1$, we can explicit the critical behavior of the free energy:
     \begin{equation}\label{critfe}
      F(\gb ) \sim C_{1} (\gb - \gb_{c}^{1})^{2},{ } \gb \searrow \gb_{c}^{1}
     \end{equation}  
      where $C_{1} :=  \frac{5}{(\sqrt{5}-1)^{2}\pi c_{K}^{2}}$. 
       \end{theorem}
       
          It is an interesting phenomenon that these explicit critical 
      points satisfy the strict (intuition matching) inequality  $\gb_{c}^{0} > \gb_{c}^{1} > \gb_{c}^{2}$.
      
       Also, the critical behavior of the free energy \eqref{critfe} matches with the one of the standard
       homogeneous wetting model (see \cite[Chapter 1, (1.20) or Chapter 2, Theorem 3]{GB}), which is given by
      \begin{equation}
       F(\gb) \sim \frac{1}{\pi c_{K}^{2}}(\gb -\gb_{c}^{0})^{2}, 
      \end{equation}  
       up to the factor $\frac{5}{(\sqrt{5}-1)^{2}} $.

      \textit{Proof}

      A byproduct of our characterization of
      the critical point of Section\ref{MRT} (see in particular Proposition \ref{ExFree}) is the fact that in this case, $e^{-\gb_{c}^{a}}$ is
      equal to the spectral
      radius of the  matrix $M_{a}$
    which is a tridiagonal $(a+1) \times (a+1)$ matrix defined by 
     \begin{equation}
      M_{a} := \begin{bmatrix} q & p & 0 & 0 & \ldots  \\ 
          p & q & p & 0 & \ldots \\
          0 & \ddots & \ddots & \ddots & 0 \\ 
          \vdots  & 0& p& q & p \\
          & \ldots& 0 & p & \frac{1+q}{2}  
               \end{bmatrix}.
     \end{equation} 
      
      When $a = 1$, this characterization directly
      leads to the equality \eqref{pq1}. 
       
      To show the critical behavior \eqref{critfe}, we first notice that the free energy $F(\gb)$ is such that 
     the spectral radius of the matrix
     \begin{equation}
      \begin{bmatrix}  q e^{-F(\gb)}  & p e^{-F(\gb)}  \\
          p e^{-F(\gb)} & \sum_{j \geq 1} K(j) e^{-F(\gb)j}  
               \end{bmatrix}
     \end{equation} 
      is equal to $e^{-(\gb - \gb_{c}^{1})}$. Proceeding as in \cite[Chapter 1]{GB} (see in particular his relations (1.16) and (1.20)), 
      we get that, as $b \searrow 0$, 
      \begin{equation}
       \frac{1+q}{2} - \sum_{j \geq 1} K(j) e^{-bj} \sim 2 c_{K} \sqrt{\pi b},
      \end{equation}  
        which we can use to deduce (since $F(\gb) \to 0$ as $\gb \to \gb_{c}^{1}$)
       \begin{equation}
        c_{K}\left(\frac{1-q}{2\sqrt{\Delta(q)}}-1\right) \sqrt{\pi F(\gb )} \sim -(\gb - \gb_{c}^{1})
       \end{equation}  
       as $\gb \searrow \gb_{c}^{1}$, where $\Delta(q) = \left(\frac{1+3q}{2} \right)^{2} - (q^{2}+4q-1) = \frac{5}{4}(q-1)^{2}$
       is the discriminant of the characteristic polynomial of the matrix $M_{1}$. 
        From this \eqref{critfe} follows.

      When $a=2$, we have to  compute the roots of the characteristic polynomial of $M_{2}$ which is given by 
      \begin{equation}
       P_{q}(X) = X^{3} - \frac{1+5q}{2}X^{2} + \frac{3q^{2}+4q-1}{2}X-\frac{q^{3}+9q^{2}-q-1}{8}. 
      \end{equation}   
     
      This is performed using Cardan-Tartaglia's formulas; one then deduces \eqref{pq2} for example with the help of a software for formal calculations. 
     $ \qed $

%
%
  
  \textbf{ Acknowledgement}: Most of the results shown here were obtained during the PhD 
   thesis of the author, and he is very grateful to Giambattista Giacomin and Francesco Caravenna for several enlightening 
  discussions about this paper. He is also very grateful for the comments of anonymous referees, who lead to a substantial improvement 
   of a previous version of this paper.

 \bibliographystyle{alpha}

  \bibliography{SPA14}

\begin{thebibliography}{CKMV09}

\bibitem[Asm03]{As}
S.~Asmussen.
\newblock {\em Applied probability and queues}, volume~51 of {\em Applications
  of Mathematics (New York)}.
\newblock Springer-Verlag, New York, second edition, 2003.
\newblock Stochastic Modelling and Applied Probability.

\bibitem[CC13]{CaCha}
F.~Caravenna and L.~Chaumont.
\newblock An invariance principle for random walk bridges conditioned to stay
  positive.
\newblock {\em Electron. J. Probab.}, 18:1--32, 2013.

\bibitem[CD08]{CD1}
F.~Caravenna and J.-D. Deuschel.
\newblock Pinning and wetting transition for {$(1+1)$}-dimensional fields with
  {L}aplacian interaction.
\newblock {\em Ann. Probab.}, 36(6):2388--2433, 2008.

\bibitem[CD09]{CD2}
F.~Caravenna and J.-D. Deuschel.
\newblock Scaling limits of {$(1+1)$}-dimensional pinning models with
  {L}aplacian interaction.
\newblock {\em Ann. Probab.}, 37(3):903--945, 2009.

\bibitem[CGZ06]{CGZ}
F.~Caravenna, G.~Giacomin, and L.~Zambotti.
\newblock Sharp asymptotic behavior for wetting models in {$(1+1)$}-dimension.
\newblock {\em Electron. J. Probab.}, 11:no. 14, 345--362 (electronic), 2006.

\bibitem[CGZ07]{CGZ1}
F.~Caravenna, G.~Giacomin, and L.~Zambotti.
\newblock Infinite volume limits of polymer chains with periodic charges.
\newblock {\em Markov Process. Related Fields}, 13(4):697--730, 2007.

\bibitem[CKMV09]{CrKMV}
M.~Cranston, L.~Koralov, S.~Molchanov, and B.~Vainberg.
\newblock Continuous model for homopolymers.
\newblock {\em J. Funct. Anal.}, 256(8):2656--2696, 2009.

\bibitem[DGZ05]{DGZ}
J.-D. Deuschel, G.~Giacomin, and L.~Zambotti.
\newblock Scaling limits of equilibrium wetting models in {$(1+1)$}-dimension.
\newblock {\em Probab. Theory Related Fields}, 132(4):471--500, 2005.

\bibitem[DIM77]{DIM}
R.~T. Durrett, D.~L. Iglehart, and D.~R. Miller.
\newblock Weak convergence to {B}rownian meander and {B}rownian excursion.
\newblock {\em Ann. Probability}, 5(1):117--129, 1977.

\bibitem[Don10]{Do3}
R.A. Doney.
\newblock Local behavior of first passage probabilities.
\newblock {\em Probab. Theory Related Fields}, 5:299--315, 2010.

\bibitem[Fun08]{FuN}
T.~Funaki.
\newblock A scaling limit for weakly pinned {G}aussian random walks.
\newblock In {\em Proceedings of {RIMS} {W}orkshop on {S}tochastic {A}nalysis
  and {A}pplications}, RIMS K\^oky\^uroku Bessatsu, B6, pages 97--109. Res.
  Inst. Math. Sci. (RIMS), Kyoto, 2008.

\bibitem[Gia07]{GB}
G.~Giacomin.
\newblock {\em Random polymer models}.
\newblock Imperial College Press, London, 2007.

\bibitem[IY01]{IY}
Y.~Isozaki and N.~Yoshida.
\newblock Weakly pinned random walk on the wall: pathwise descriptions of the
  phase transition.
\newblock {\em Stochastic Process. Appl.}, 96(2):261--284, 2001.

\bibitem[MP10]{MP}
P.~M{\"o}rters and Y.~Peres.
\newblock {\em Brownian motion}.
\newblock Cambridge Series in Statistical and Probabilistic Mathematics.
  Cambridge University Press, Cambridge, 2010.
\newblock With an appendix by Oded Schramm and Wendelin Werner.

\bibitem[MT09]{MT}
S.~Meyn and R.~L. Tweedie.
\newblock {\em Markov chains and stochastic stability}.
\newblock Cambridge University Press, Cambridge, second edition, 2009.
\newblock With a prologue by P. W. Glynn.

\bibitem[Shi83]{Shi}
M.~Shimura.
\newblock A class of conditional limit theorems related to ruin problem.
\newblock {\em Ann. Probab.}, 11(1):40--45, 1983.

\bibitem[Soh13]{Soh1}
J.~Sohier.
\newblock The scaling limits of a heavy tailed markov renewal process.
\newblock {\em Ann. Inst. H. Poincar\'e Probab. Statist.}, 2:483--505, 2013.

\bibitem[Zer87]{Ze}
M.~Zerner.
\newblock Quelques propri\'et\'es spectrales des op\'erateurs positifs.
\newblock {\em J. Funct. Anal.}, 72(2):381--417, 1987.

\end{thebibliography}

  \end {document}